\documentclass[12pt]{article}
\usepackage{amssymb,epsfig}

\usepackage[french]{babel}
\usepackage[T1]{fontenc}
\usepackage{color}

\newtheorem{defi}{D\'efinition}[section]
\newtheorem{prop}[defi]{Proposition}
\newtheorem{theo}[defi]{Th\'eor\`eme}
\newtheorem{conj}[defi]{Conjecture}
\newtheorem{lemm}[defi]{Lemme}
\newtheorem{coro}[defi]{Corollaire}
\newtheorem{rema}[defi]{Remarque}
\newtheorem{exem}[defi]{Exemple}
\newtheorem{exer}{Exercice}

\newcommand{\bdefi}{\begin{defi}}
\newcommand{\edefi}{\end{defi}}
\newcommand{\bprop}{\begin{prop}}
\newcommand{\eprop}{\end{prop}}
\newcommand{\btheo}{\begin{theo}}
\newcommand{\etheo}{\end{theo}}
\newcommand{\blemm}{\begin{lemm}}
\newcommand{\brema}{\begin{rema}}
\newcommand{\erema}{\end{rema}}
\newcommand{\bexer}{\begin{exer}}
\newcommand{\eexer}{\end{exer}}
\newcommand{\bconj}{\begin{conj}}
\newcommand{\econj}{\end{conj}}
\newcommand{\elemm}{\end{lemm}}
\newcommand{\bcoro}{\begin{coro}}
\newcommand{\ecoro}{\end{coro}}
\newcommand{\bexem}{\begin{exem}}
\newcommand{\eexem}{\end{exem}}
\newcommand{\dem}{\noindent{\bf Preuve. }}

\newcommand{\B}{{\cal B}}

\newcommand{\G}{{\cal G}}

\newcommand{\D}{{\cal D}}

\renewcommand{\H}{{\cal H}}
\renewcommand{\O}{{\cal O}}

\usepackage{amssymb}

\newcommand{\maths}[1]{{\mathbb #1}}  

\newcommand{\RR}{\maths{R}}
\newcommand{\NN}{\maths{N}}

\newcommand{\QQ}{\maths{Q}}

\newcommand{\HH}{\maths{H}}

\newcommand{\FF}{\maths{F}}

\newcommand{\ZZ}{\maths{Z}}

\newcommand{\PP}{\maths{P}}

\newcommand{\TT}{\maths{T}}

\newcommand{\ra}{\rightarrow}

\newcommand{\wt}[1]{{\widetilde{#1}}}
\newcommand{\wh}[1]{{\widehat{#1}}}

\newcommand{\ga}{\gamma}
\newcommand{\Ga}{\Gamma}

\newcommand{\cqfd}{\hfill$\Box$}
\newcommand{\bac}{{\backslash\!\backslash}}

\newcommand{\vy}{v_\infty}
\newcommand{\dy}{d_\infty}
\newcommand{\va}[1]{\left|#1\right|_{\infty}}
\newcommand{\gatq}{{\Gamma\backslash\TT_q}}

\newcounter{fig}

\usepackage{epsfig}

\title{Dynamique sur le rayon modulaire \\ et \\
fractions continues en caract\'eristique $p$.}

\author{Anne Broise \and Fr\'ed\'eric Paulin}

\date{}

\begin{document}
\maketitle

\begin{center} {\bf Abstract}
\end{center}

\begin{center}\begin{minipage}[b]{13cm}{
~~~ Let $\wh K$ be the field of formal Laurent series in $X^{-1}$ over the
finite field $k$, and let $A$ be the ring of polynomials in $X$ over
$k$. One of the main results of the paper is to give a particularly
nice coding of the geodesic flow on the quotient of the Bruhat-Tits
tree $\TT$ of ${\rm PGL}_2(\wh K)$ by ${\rm PGL}_2(A)$, by using the
continued fraction expansion of the endpoints of the geodesic lines in
$\TT$ (the space of ends of $\TT$ identifies with $\PP_1(\wh K)$).
This allows in particular to prove in a dynamical way the invariance
of the Haar measure by the Artin map.  }\end{minipage}
\footnote{ {\bf AMS codes:} 11
  J 70, 20 G 25, 20 E 08, 37 A 45, 11 K 50. {\bf Keywords:} continued
  fractions, Artin map, Laurent series field, Bruhat-Tits tree,
  geodesic flow, coding.  }
\end{center}

\section{Introduction}
\label{sec:Intro}

Soit $\wh K=\FF_q((X^{-1}))$ le corps local des s\'eries formelles de
Laurent en $X^{-1}$ \`a coefficients dans le corps fini $\FF_q$, et
$\O=\FF_q[[X^{-1}]]$ le sous-anneau de $\wh K$ des s\'eries formelles
enti\`eres. On consid\`ere le groupe localement compact $G={\rm
PGL}(2,\wh K)$, et son r\'eseau non uniforme $\Ga={\rm
PGL}(2,\FF_q[X])$.  Le groupe $G$ agit sur son arbre (localement fini)
de Bruhat-Tits $\TT_q$. L'espace des bouts de $\TT_q$ s'identifie avec la
droite projective $\PP_1(\wh K)$, et on note $x_*$ le point base
standard de $\TT_q$ (voir par exemple \cite{Ser2} ou la partie
\ref{sec:notations} pour des rappels).

L'un des buts principaux de cet article, qui fait suite \`a
\cite{Pau}, est d'expliciter en termes arithm\'etiques la structure
ergodique des actions commutantes de $\Ga$ et du flot g\'eod\'esique
(action de $\ZZ$ par translation \`a la source) sur l'espace des
g\'eod\'esiques de $\TT_q$ (i.e.~des isom\'etries $\ell:\RR\ra \TT_q$
d'origine $\ell(0)$ un sommet de $\TT_q$).

Nous d\'ecrivons (dans la partie
\ref{subsec:lienartin}) la structure de l'ensemble des g\'eod\'esiques
de $\TT_q$ modulo l'action de $\Ga$.  Notons $\pi:\TT_q\ra
\Ga\backslash \TT_q$ la projection canonique. Alors, l'ensemble
$\pi^{-1}(\pi(x_*))$ est une section $\Ga$-\'equivariante globale pour
le flot g\'eod\'esique (toute orbite la rencontre une infinit\'e de
fois). De plus, toute g\'eod\'esique de $\TT_q$, d'origine dans la
section globale $\pi^{-1}(\pi(x_*))$, est \'equivalente, modulo
l'action d'un \'el\'ement de $\Gamma$, \`a une unique g\'eod\'esique
$\ell$ d'origine $x_*$, d'extr\' emit\'e n\'egative un point de
$J={\displaystyle \bigcup_{a\in\FF_q[X]-\FF_q} (a+ X^{-1}\O)}$ et
d'extr\' emit\'e positive un point de $X^{-1}\O$. Une mani\`ere de
rendre cet \'el\'ement de $\Ga$ unique est d'introduire des
``d\'ecorations'' sur les g\'eod\'esiques  (voir la
partie \ref{sec:codage}).

Notons $\G'_0(\TT_q)$ l'ensemble des telles g\'eod\'esiques $\ell$,
identifi\'ees \`a leurs couples d'extr\'emit\'es $(\xi_-,\xi_+)$, avec
de plus $\xi_+$, $\xi_-$ irrationnelles (i.e.~dans $\PP_1(\wh
K)-\PP_1(K)$).  Soit $\wt \Psi:\G'_0(\TT_q)\ra \G'_0(\TT_q)$
l'application induite par l'application de premier retour du flot
g\'eod\'esique sur la section globale $\pi^{-1}(\pi(x_*))$ (voir la
partie \ref{subsec:lienartin}).  Pour $\ell$ dans $\G'_0(\TT_q)$,
notons $(a_n)_{n\geq 1}$ le d\'eveloppement en fractions continues
d'Artin \cite{Art} de $\xi_+$, et $(a_{-n})_{n\geq 0}$ celui de
${-1\over \xi_-}$ (voir la partie \ref{subsec:corps} pour des
rappels).

Nous montrons dans la partie \ref{sec:codage} le r\'esultat suivant~:

\btheo \label{theo:intro}
L'application $\Theta':\G'_0(\TT_q)\ra (\FF_q[X]-\FF_q)^\ZZ$,
d\'efinie par $\Theta'(\ell)=(a_n)_{n\in\ZZ}$, est un
hom\'eomorphisme qui rend le diagramme suivant commutatif
$$\begin{array}[b]{ccc} \G_0'(\TT_q) &
     \stackrel{\Theta'}{\longrightarrow} & (\FF_q[X]-\FF_q)^{\ZZ} \\
     \wt\Psi\downarrow\;  & & \downarrow \sigma \\
     \G_0'(\TT_q) &\stackrel{\Theta'}{\longrightarrow} &
     (\FF_q[X]-\FF_q)^{\ZZ}
\end{array}\;,$$
o\`u $\sigma$ est le d\'ecalage \`a gauche des suites bilat\`eres de
$(\FF_q[X]-\FF_q)^{\ZZ}$.  De plus, le diagramme suivant commute
$$\begin{array}[b]{ccc}
\G'_0(\TT_q) & \stackrel{\wt\Psi}{\longrightarrow} & \G'_0(\TT_q) \\
{\rm p}_+\downarrow\;  & & \;\downarrow {\rm p}_+ \\
X^{-1}\O\cap (\wh K- K) &\stackrel{\Psi}{\longrightarrow} &
X^{-1}\O\cap (\wh K- K)
\end{array}\;,$$
o\`u $\Psi:\xi\mapsto {1\over \xi}-\left[{1\over \xi}\right]$ est
l'application d'Artin (si $\zeta$ est dans $\wh K$, alors $[\zeta]$
d\'esigne sa partie enti\`ere), et ${\rm p}_+:(\xi_-,\xi_+)\mapsto
\xi_+$ est la projection sur la deuxi\`eme coordonn\'ee.
\etheo

Notons $m$ la mesure de probabilit\'e qui est la restriction \`a
$\G_0'(\TT_q)$ de la mesure naturelle invariante par le flot
g\'eod\'esique. La mesure $m$ co\"{\i}ncide avec la mesure de Haar
quotient sur $\Ga\backslash G/M$ pour $M$ le sous-groupe des matrices
diagonales \`a coefficients dans $\O^*$, et, de manière classique,
avec la mesure de Patterson-Sullivan-Bowen-Margulis de $\Ga$ (voir par
exemple \cite{Bou, BM} et la partie \ref{sec:mesure}).

Nous montrons dans la partie \ref{sec:mesure} que l'image de $m$ par
$\Theta'$ est une mesure de Bernoulli sur $(\FF_q[X]-\FF_q)^{\ZZ}$.
Ceci implique en particulier que si $S$ est le sous-groupe diagonal de
$G$, alors l'action \`a droite de $S/M$ sur $\Ga\backslash G/M$ est
Bernoulli (donc m\'elangeante, ce qui \'etait d\'ej\`a connu, voir par
exemple \cite{BN}). Dans un article en pr\'eparation \cite{BP}, nous
\'etudierons le cas g\'en\'eral des r\'eseaux des groupes
alg\'ebriques semi-simples de rang $1$ sur un corps local non
archim\'edien. En fait, nous donnerons des codages markoviens de flots
g\'eod\'esiques sur des arbres munis d'actions tr\`es g\'en\'erales de
groupes.  Ces codages permettent de contourner l'abondance de torsion
dans les r\'eseaux non-uniformes d'arbres, dont on ne peut se
  débarrasser par passage à un sous-groupe d'indice fini.

Nous montrons dans la partie \ref{sec:applications} que l'image de $m$
par la deuxi\`eme projection est la mesure de Haar sur
$X^{-1}\O$. Ceci explique de mani\`ere dynamique l'invariance de cette
mesure de Haar par la transformation d'Artin.

Ces r\'esultats sont analogues aux r\'esultats qui relient le flot
g\'eod\'esique sur la courbe mo\-du\-laire ${\rm PSL}(2,\ZZ)\backslash
\HH^2$ avec le d\'eveloppement en fractions continues des nombres
r\'eels (et qui expliquent, en particulier, l'invariance de la mesure
de Gauss $\frac{1}{\log 2}\frac{dx}{1+x}$ par la transformation de
Gauss $x\mapsto {1\over x}-[{1\over x}]$) (voir par exemple
\cite{Seri}).

Cet article fait suite \`a \cite{Pau}, o\`u une partie de l'analogie
ayant trait aux g\'eod\'esiques individuelles, est d\'evelopp\'ee.
Mais le codage global du flot g\'eod\'esique n'est pas contenu (m\^eme
pas entre les lignes) dans \cite{Pau}, car une approche globale
suivant de trop pr\`es \cite{Pau} conduit \`a des discontinuit\'es. La
pr\'esence du corps r\'esiduel $\FF_q$, absent dans le cas r\'eel, est
une des sources de probl\`emes. C'est, entre autre, le travail de
``naturalit\'e'' du pr\'esent article, en particulier \`a partir de la
notion de g\'eod\'esique d\'ecor\'ee, qui permet le codage global. Il
ne s'agit pas de construire n'importe quel codage, mais un qui soit
intimement lié à la structure arithmétique du réseau ${\rm
  PGL}(2,\FF_q[X])$ (celui-ci n'est pas de type fini et contient des
sous-groupes infini de torsion), et qui permette une correspondance
entre propriétés dynamiques et arithmétiques.

\section{Notations et rappels}
\label{sec:notations}

Toute cette partie est compos\'e de rappels, pour lesquels nous
renvoyons par exemple \`a \cite{Ser1,Spr,Sch,Laj,BN,Ser2,Pau}.  Elle
n'est \'ecrite que pour \'eviter au lecteur qui ne conna\^{\i}trait
pas les notations et r\'esultats de \cite{Ser1,Pau} d'avoir \`a lire
le pr\'esent article en ayant \`a c\^ot\'e ces deux
r\'ef\'erences. Les autres lecteurs peuvent se reporter directement au
chapitre \ref{sec:codage}.

\subsection{Le corps des s\'eries formelles de Laurent}
\label{subsec:corps}

Soit $k=\FF_q$ un corps fini, d'ordre $q$ (o\`u $q$ est une
puissance d'un nombre premier). On note $A=k[X]$ l'anneau des
polyn\^omes en une variable $X$ sur $k$, et $K=k(X)$ le corps des
fractions rationnelles en $X$ sur $k$.  Soit $\wh K=k((X^{-1}))$ le
compl\'et\'e de $K$ pour la valuation $\vy$ d\'efinie par
$$\vy(P/Q)=(-\deg P)-(-\deg Q)\;.$$
Le corps $\wh K$ est muni de l'unique valuation qui \'etend $\vy$ (que
l'on notera de la m\^eme mani\`ere), de la valeur absolue
$$\va{f}=q^{-\vy(f)} $$
et de la distance ultram\'etrique d\'efinie
par cette valeur absolue
$$\dy(f,g)=\va{f-g}.$$
Tout \'el\'ement de $\wh K$ s'\'ecrit de mani\`ere unique
comme s\'erie convergente
$$f=\sum_{i>>-\infty} f_iX^{-i}$$
avec $f_i$ dans $k$, nul pour $i$ suffisamment petit. On a
$$\vy(f)=\sup\{j\in\ZZ\;|\; \;\forall i<j, f_i=0\}\;.$$

On note $\O=\{f\in\wh K, \vy(f)\geq 0\}$ l'anneau de la valuation
$\vy$ dans $\wh K$. C'est le sous-espace compact-ouvert $k[[X^{-1}]]$
des s\'eries enti\`eres en $X^{-1}$ sur $k$, c'est aussi la boule
ferm\'ee de rayon $1$ et de centre $0$ dans $\wh K$.

Pour tout $f$ dans $\wh K$, il existe un unique couple form\'e d'un
polyn\^ome en $X$ sur $k$, not\'e $[f]\in A$, et d'une s\'erie
enti\`ere en $X^{-1}$ sur $k$, de terme constant nul, not\'ee
$\{f\}\in X^{-1}\O$, tels que $$f=[f]+\{f\}.$$
On appelle $[f]$ la {\it partie enti\`ere} de $f$ et $\{f\}=f-[f]$ la
{\it partie fractionnaire} de $f$.  {\it L'application d'Artin} est
l'application $\Psi:X^{-1}\O-\{0\}\ra X^{-1}\O$ d\'efinie par
$$\Psi(f)=\{\frac{1}{f}\}=\frac{1}{f}-[\frac{1}{f}].$$

On note $^c\! K$ l'ensemble $\wh K-K$ des points irrationnels de $\wh
K$. Pour $f$ dans $^c\! K$, on pose $a_0=[f]$ et pour $n\geq 1$
entier,
$$a_n=\left[\frac{1}{\Psi^{n-1}(f-a_0)}\right].$$
Alors $a_n $ est dans $A$. Si $n\geq 1$, le degr\'e de $a_n$ est
strictement positif et
$$f=\lim_{n\ra +\infty}
a_0 +\frac{1}{\displaystyle a_1+\frac{1}{\;\;
\begin{array}{cccc}\cdot & & & \\
& \cdot & &  \\
&  &  \cdot & \\
&  &  & a_{n-1} +\frac{1}{a_n}
\end{array}}}.$$

Pour $f$ dans $X^{-1}\O\cap\,^c\! K$, on a $a_0=0$ et on appelle {\it
     d\'eveloppement en fractions continues d'Artin} de $f$ la suite
$(a_n)_{n\geq 1}$.

\subsection{L'arbre de Bruhat-Tits}
\label{subsec:arbre}

Notons $G$ le groupe localement compact ${\rm PGL}(2,\wh{K})$. On
appelle {\it groupe modulaire} le sous-groupe discret ${\rm PGL}(2,A)$
de $G$, et on le note $\Gamma$.  Dans toute la suite, on note de la
m\^eme mani\`ere une matrice de ${\rm GL}(2,\wh{K})$ et son image dans
$G$.

On rappelle que pour tout corps commutatif $\kappa$, l'action par
homographies du groupe ${\rm PGL}(2,\kappa)$ sur la droite projective
$\PP_1(\kappa)$ est simplement transitive sur les triplets de points
de $\PP_1(\kappa)$. Comme l'application naturelle ${\rm PGL}(2,A)\ra{\rm
    PGL}(2,K)$ est une bijection, le groupe $\Ga$ agit simplement
transitivement sur les triplets de points de $\PP_1(K)$.

On note ${\rm GL}(2,\wh{K})_1$ le groupe des matrices carr\'ees de
taille $2$ \`a coefficients dans $\wh K$, dont la valeur absolue du
d\'eterminant est \'egale \`a $1$, et $G_1={\rm PGL}(2,\wh{K})_1$ le
quotient par son centre.  Comme les \'el\'ements inversibles de $A$
sont de valeur absolue \'egale \`a $1$, le groupe $\Gamma$ est contenu
dans $G_1$.

\medskip
L'{\it arbre de Bruhat-Tits} $\TT_q$ de $(SL_2,\wh K)$ est le graphe
d\'efini par
\begin{enumerate}
\item les sommets de $\TT_q$ sont les classes d'homoth\'etie (par $\wh
     K^\times$) $\Lambda=[L]$ de $\O$-r\'eseaux (i.e.~$\O$-sous-modules
     libres de rang deux) $L$ dans $\wh K\times\wh K$;
\item deux sommets $\Lambda,\Lambda'$ sont joints par une ar\^ete si
     et seulement s'il existe des re\-pr\'e\-sen\-tants $L,L'$ de
     $\Lambda,\Lambda'$ tels que $L'\subset L$ et $L/L'$ est isomorphe
     \`a $\O/X^{-1}\O=k$.
\end{enumerate}

On note $V\TT_q$ l'ensemble des sommets de $\TT_q$ et $E\TT_q$
l'ensemble de ses ar\^etes. On note $d$ la distance dans $\TT_q$ et
$x_*=[\O\times\O]$ la classe du r\'eseau standard. Le graphe $\TT_q$
est un arbre r\'egulier de degr\'e $q+1$. L'action naturelle de ${\rm
    GL}(2,\wh K)$ sur les $\O$-r\'eseaux induit une action de $G$ sans
inversion sur $\TT_q$, transitive sur les ar\^etes.  Le stabilisateur
de $x_*$ est le groupe ${\rm PGL}(2,\O)$.

On identifie dans la suite l'arbre $\TT_q$ avec sa r\'ealisation
g\'eom\'etrique. On note $\partial \TT_q$ l'espace des bouts de
l'espace topologique localement compact $\TT_q$. C'est l'espace des
classes d'\'equivalence de rayons g\'eod\'esiques dans $\TT_q$, o\`u
l'on identifie deux rayons g\'eod\'esiques si leur intersection est
encore un rayon g\'eod\'esique. L'action de $G$ sur $\TT_q$ s'\'etend
continuement en une action par hom\'eomorphismes de $G$ sur $\partial
\TT_q$.

L'espace $\partial \TT_q$ s'identifie, de mani\`ere
$G$-\'equivariante, avec la droite projective $\PP(\wh K\times\wh
K)=\PP_1(\wh K)$, par l'application qui, \`a l'extr\'emit\'e d'un
rayon g\'eod\'esique issu de $\Lambda_0=x_*$, de suite des sommets
cons\'ecutifs $(\Lambda_n)_{n\in\NN}$, associe l'unique droite de $\wh
K\times\wh K$ contenant l'intersection des $\O$-r\'eseaux $L_n$, avec
$(L_n)_{n\in\NN}$ l'unique suite de $\O$-r\'eseaux dans $\wh
K\times\wh K$ telle que $[L_n]=\Lambda_n$ et $L_{n+1}\subset L_n$.

En choisissant la droite $\wh K\times\{0\}$ de $\wh K\times\wh K$
comme point \`a l'infini de $\PP_1(\wh K)$, on identifie $\partial
\TT_q=\PP_1(\wh K)$ avec $\wh K\cup\{\infty\}$. L'action de $G$ sur
$\partial \TT_q=\PP_1(\wh K)$ correspond \`a l'action de $G$ par
homographies sur $\wh K\cup\{\infty\}$.

\medskip
On note $\pi:\TT_q\ra \Ga\backslash\TT_q$ la projection
canonique et $\Ga_\infty$ le fixateur dans $\Ga$ du point $\infty$ de
$\partial \TT_q$. Le rayon g\'eod\'esique $\D_\Gamma$ issu de $x_*$ et
d'extr\'emit\'e $\infty$, donc de suite des sommets cons\'ecutifs
$([\O\times X^{-n}\O])_{n\in\NN}$, est un domaine fondamental pour
l'action du groupe modulaire $\Gamma$ sur l'arbre de Bruhat-Tits
$\TT_q$, au sens o\`u ses images par $\Gamma$ recouvrent $\TT_q$.  La
restriction $\pi\!\mid_{\D_\Gamma}$ est un isomorphisme simplicial de
$\D_\Gamma$ sur $\Ga\backslash\TT_q$. Le graphe $\Ga\backslash\TT_q$
h\'erite par $\pi$ d'une structure de graphe de groupes, que l'on note
$\Gamma\bac\TT_q$, et que l'on appelle le {\it rayon modulaire} (voir
\cite{Ser2} 
pour les d\'efinitions et la construction).

Voir par exemple \cite{Ser2} pour des justifications et compl\'ements.

\subsection{La famille $\Gamma$-\'equivariante maximale d'horoboules
     d'int\'erieurs disjoints.}
\label{subsec:horo}

Pour tout $\xi$ de $\partial \TT_q$, on appelle {\it fonction de
    Buseman} de $\TT_q$ l'application $\beta_\xi:\TT_q\times
\TT_q\ra\RR$ d\'efinie par
$$\beta_\xi(x,y)=\lim_{t\ra +\infty} d(y,c(t))-d(x,c(t))$$
o\`u $c:[0,+\infty[\,\ra \TT_q$ est un rayon g\'eod\'esique
convergeant vers $\xi$.  La fonction $\beta_\xi$ ne d\'epend pas du
choix de $c$, est invariante par isom\'etries~:
$$\forall \ga\in G, \;\;\;\;\beta_{\ga\xi}(\ga x,\ga
y)=\beta_\xi(x,y),$$
et v\'erifie la relation de cocycle~:
$$\beta_\xi(x,y)+\beta_\xi(y,z)=\beta_\xi(x,z).$$
En fait, l'application $t\mapsto d(y,c(t))-d(x,c(t))$ est constante
\`a partir d'un certain temps. Le rayon g\'eod\'esique $[x,\xi[$
rencontre le rayon g\'eod\'esique $[y,\xi[$ en un rayon g\'eod\'esique
$[w,\xi[$ et $\beta_\xi(x,y)=d(y,w)-d(w,x)$.

On appelle {\it horosph\`ere} centr\'ee en $\xi\in\partial \TT_q$ et
passant par $x\in\TT_q$ l'ensemble des points $y$ de $ \TT_q$ tels que
$\beta_\xi(x,y)=0$. Son {\it horoboule} associ\'ee est l'ensemble des
points $y$ de $ \TT_q$ tels que $\beta_\xi(x,y)\leq 0$.

\begin{center}
\begin{picture}(0,0)%
\includegraphics{fig_famillehoro.pstex}%
\end{picture}%
\setlength{\unitlength}{3947sp}%
\begingroup\makeatletter\ifx\SetFigFont\undefined%
\gdef\SetFigFont#1#2#3#4#5{%
  \reset@font\fontsize{#1}{#2pt}%
  \fontfamily{#3}\fontseries{#4}\fontshape{#5}%
  \selectfont}%
\fi\endgroup%
\begin{picture}(4459,3577)(156,-3148)
\put(1051,-1336){\makebox(0,0)[lb]{\smash{\SetFigFont{12}{14.4}{\rmdefault}{\mddefault}{\updefault}{$\infty$}%
}}}
\put(1966,-2771){\makebox(0,0)[lb]{\smash{\SetFigFont{12}{14.4}{\rmdefault}{\mddefault}{\updefault}{$\ga\infty$}%
}}}
\put(4246,-1346){\makebox(0,0)[lb]{\smash{\SetFigFont{12}{14.4}{\rmdefault}{\mddefault}{\updefault}{$0$}%
}}}
\put(3271,-1976){\makebox(0,0)[lb]{\smash{\SetFigFont{12}{14.4}{\rmdefault}{\mddefault}{\updefault}{$H_0$}%
}}}
\put(2226,-286){\makebox(0,0)[lb]{\smash{\SetFigFont{12}{14.4}{\rmdefault}{\mddefault}{\updefault}{$H_1$}%
}}}
\put(2451,-1986){\makebox(0,0)[lb]{\smash{\SetFigFont{12}{14.4}{\rmdefault}{\mddefault}{\updefault}{$H_{\ga\infty}$}%
}}}
\put(4551,-576){\makebox(0,0)[lb]{\smash{\SetFigFont{12}{14.4}{\rmdefault}{\mddefault}{\updefault}{$X^{-1}\O$}%
}}}
\put(2656,-1486){\makebox(0,0)[lb]{\smash{\SetFigFont{12}{14.4}{\rmdefault}{\mddefault}{\updefault}{$x_*$}%
}}}
\put(1806,-921){\makebox(0,0)[lb]{\smash{\SetFigFont{12}{14.4}{\rmdefault}{\mddefault}{\updefault}{$H_\infty$}%
}}}
\put(2681,264){\makebox(0,0)[lb]{\smash{\SetFigFont{12}{14.4}{\rmdefault}{\mddefault}{\updefault}{$1$}%
}}}
\put(711,-2466){\makebox(0,0)[lb]{\smash{\SetFigFont{12}{14.4}{\rmdefault}{\mddefault}{\updefault}{$J=$}%
}}}
\put(2816,-956){\makebox(0,0)[lb]{\smash{\SetFigFont{12}{14.4}{\rmdefault}{\mddefault}{\updefault}{$v_*$}%
}}}
\put(156,-2811){\makebox(0,0)[lb]{\smash{\SetFigFont{12}{14.4}{\rmdefault}{\mddefault}{\updefault}{$\displaystyle\bigcup_{a\in A-k} a+X^{-1}\O$}%
}}}
\end{picture}

\end{center}
\begin{center}
    {\bf Figure \addtocounter{fig}{1}\arabic{fig} :} Trajet des
    g\'eod\'esiques dans la famille d'horoboules $\H\B$.
\end{center}

On note $H_\infty$ l'horosph\`ere centr\'ee en $\infty$ de
$\partial\TT_q$ et passant par $x_*$, et $H\!B_\infty$ son horoboule
associ\'ee.  C'est l'orbite du rayon fondamental $\D_\Gamma$ par le
fixateur $\Gamma_\infty$ du point $\infty$ (de $\partial \TT_q$) dans
le groupe modulaire $\Gamma$. Pour tout $\gamma$ dans $\Gamma$, on
note $H_{\gamma\infty} =\gamma H_\infty$, c'est l'horosph\`ere
centr\'ee en $\gamma\infty$ et passant par $\gamma x_*$, et on note
$H\!B_{\gamma\infty} =\gamma H\!B_\infty$ son horoboule associ\'ee.
Il est facile de voir que la famille
$\H\B=(H\!B_{\gamma\infty})_{\gamma\in\Gamma/\Gamma_\infty}$ est la
famille des adh\'erences des composantes connexes de la pr\'eimage par
$\pi$ du rayon $\gatq$ priv\'e de son origine. Les horoboules de cette
famille se rencontrent deux \`a deux en au plus un point (sur
l'intersection des horosph\`eres associ\'ees), et leur r\'eunion est
\'egale \`a $\TT_q$.

Notons $\H$ la famille d'horosph\`eres
$(H_{\gamma\infty})_{\gamma\in\Gamma/\Gamma_\infty}$.  (Voir par
exemple \cite{Pau,Pau2} pour des justifications et compl\'ements.)

\section{Codage du flot g\'eod\'esique sur le rayon modulaire}
\label{sec:codage}

\subsection{L'espace des g\'eod\'esiques d\'ecor\'ees
sur le rayon modulaire}
\label{subsec:geoddeco}

Appelons {\it g\'eod\'esique d\'ecor\'ee de $\TT_q$} une
g\'eod\'esique $\ell$ de $\TT_q$ (i.e.~une isom\'etrie
$\ell:\RR\ra\TT_q$), dont l'origine $x_0=x_0(\ell)=\ell(0)$ est sur
l'une des horosph\`eres de la famille $\H$ (i.e.~est un point de $\pi
^{-1}(\pi (x_*))$ ), et qui est munie d'une ar\^ete $v_0=v_0(\ell)$
issue de $x_0$ non contenue dans $\ell$, appel\'ee {\it d\'ecoration}.
Remarquons que pour $q=2$, la d\'ecoration est unique.

La {\it g\'eod\'esique d\'ecor\'ee standard} $\ell_*$ est la
g\'eod\'esique d'extr\'emit\'es $\infty$ et $0$, orient\'ee de
$\infty$ vers $0$, d'origine $x_*$ et d\'ecor\'ee par l'ar\^ete $v_*$
qui pointe vers 1.

Notons $\G_0(\TT_q)$ l'ensemble des g\'eod\'esiques d\'ecor\'ees {\it
    totalement irrationnelles} de $\TT_q$, c'est-\`a-dire celles qui ont
leurs deux extr\'emit\'es dans $\partial\TT_q-\PP_1(K)=\,^c\!K$.
Attention, $\ell_*$ est une g\'eod\'esique d\'ecor\'ee qui n'est pas
dans $\G_0(\TT_q)$.


L'ensemble $\G_0(\TT_q)$ est muni de la topologie d\'efinie par le
syst\`eme fondamental d'en\-tou\-ra\-ges $\{V_n\}_{n\in \NN}$, o\`u $V_n$
est l'ensemble des couples de g\'eod\'esiques d\'ecor\'ees qui
co{\"\i}ncident entre les instants $-n$ et $n$ et qui ont m\^eme
d\'ecoration.

\blemm
Le groupe $\Ga$ agit librement et proprement sur $\G_0(\TT_q)$.
\elemm

\dem
Le groupe $\Ga$ pr\'eserve la famille $\H$ et
$\PP_1(K)$ donc il pr\'eserve $\,^c\!K$. Il agit sur $\G_0(\TT_q)$~: l'image
par un \'el\'ement $\ga$ de $\Ga$ d'une g\'eod\'esique $\ell$ de
d\'ecoration $v_0$ est la g\'eod\'esique $\ga\ell$ de d\'ecoration
$\ga v_0$.  Comme $\Ga$ agit proprement sur $\TT_q$, l'action de $\Ga$
sur $\G_0(\TT_q)$ est propre. Cette action est aussi libre, car $G$
agit simplement transitivement sur les triplets de points de $\partial
\TT_q$, et si un \'el\'ement $\ga$ de $\Ga$ fixe une g\'eod\'esique
d\'ecor\'ee, d'extr\'emit\'es (distinctes) $\xi_-,\xi_+$ et de
d\'ecoration $v_0$, alors elle fixe aussi le point \`a l'infini
$\xi_0$ de l'horoboule de la famille $\H\B$ qui contient $v_0$, et
$\xi_0$ est distinct de $\xi_-$ et de $\xi_+$.
\cqfd

\medskip
Notons $\G_0(\Ga\bac\TT_q)$ le quotient $\Ga\backslash\G_0(\TT_q)$.

\bigskip Pour $\ell$ une g\'eod\'esique de $\TT_q$, on note
$\xi_-=\xi_-(\ell)$ et $\xi_+=\xi_+(\ell)$ ses extr\'emit\'es
n\'egative et positive. On munit $\partial\TT_q$ de la topologie
d\'efinie par $d_\infty$, $V\TT_q$ et $E\TT_q$ de la topologie
discr\`ete et $\partial\TT_q\times \partial\TT_q\times V\TT_q\times
E\TT_q$ de la topologie produit.  L'application qui \`a un \'el\'ement
$\ell$ de $\G_0(\TT_q)$ associe le quadruplet $(\xi_-,\xi_+,x_0,v_0)$
de $\partial\TT_q\times \partial\TT_q\times V\TT_q\times E\TT_q$ est
alors un hom\'eomorphisme de $\G_0(\TT_q)$ sur son image, qui est un
bor\'elien $B$. Nous appellerons cette application le {\it
    param\'etrage} de l'espace des g\'eod\'esiques d\'ecor\'ees
totalement irrationnelles de $\TT_q$. Par la suite, nous identifierons
une g\'eod\'esique d\'ecor\'ee $\ell$ et son quadruplet
$(\xi_-,\xi_+,x_0,v_0)$. L'action de $\Gamma$ sur $\G_0(\TT_q)$
s'identifie alors avec la restriction \`a $B$ de l'action diagonale de
$\Gamma$ sur l'espace produit $\partial\TT_q\times \partial\TT_q\times
V\TT_q\times E\TT_q$.

\subsection{Codage du flot g\'eod\'esique sur le rayon modulaire}
\label{sec:codageun}

Appelons {\it renversement du temps sur $\G_0(\TT_q)$} l'application
$\wt\tau:\G_0(\TT_q)\to \G_0(\TT_q)$ d\'efinie par
$\wt\tau(\xi_-,\xi_+,x_0,v_0)=(\xi_+,\xi_-,x_0,v_0)$.  Il est clair
que $\wt\tau$ est \'equivariante sous l'action de $\Gamma$, elle
induit donc par passage au quotient une application
$\tau:\G_0(\Gamma\bac\TT_q)\to \G_0(\Gamma\bac\TT_q)$, que nous
appelerons {\it renversement du temps sur $\G_0(\Gamma\bac\TT_q)$}.

Nous allons maintenant d\'efinir une application $T:\G_0(\Gamma\bac\TT_q)
\ra\G_0(\Gamma\bac\TT_q)$, que l'on peut voir comme
l'application de premier retour en l'origine du flot g\'eod\'esique
sur le rayon modulaire.

Soit $\ell=(\xi_-,\xi_+,x_0,v_0)$ dans $\G_0(\TT_q)$. Notons
$\xi_n=\xi_n(\ell)$ les points de $\PP_1(K)$ tels que
$(H\!B_{\xi_n})_{n\in\ZZ}$ soit la suite des horoboules de $\H\B$
travers\'ees cons\'ecutivement par $\ell$, avec $H\!B_{\xi_0}\cap
H\!B_{\xi_1}=\{x_0\}$. Soit $x_n=x_n(\ell)$ le point de $V\TT_q$ tel
que $H\!B_{\xi_n}\cap H\!B_{\xi_{n+1}}=\{x_n\}$.

Il faut remarquer que les suites $(\xi_n)$ et $(x_n)$ sont bien
index\'ees par $\ZZ$. En effet, par d\'efinition, une g\'eod\'esique
d\'ecor\'ee de $\G_0(\TT_q)$ a ses deux extr\'emit\'es irrationnelles.
Or toute g\'eod\'esique qui rentre dans une horoboule de la famille
$\H\B$ ou bien en ressort, ou bien converge vers le point \`a l'infini
de cette horoboule, qui est rationnel.

Notons $\eta_0=\eta_0(\ell)$ le point de $\PP_1(K)$ tel que
l'horoboule $H\!B_{\eta_0}$ contienne l'ar\^ete $v_0$. Soit
$s_0=s_0(\ell)$ l'unique \'el\'ement (d'ordre $2$) de $\Gamma$ qui fixe
$\eta_0$ et \'echange $\xi_{0}$ et $\xi_1$. Soit
$s_{\frac{1}{2}}=s_{\frac{1}{2}}(\ell)$ l'unique \'el\'ement (d'ordre
$2$) de $\Gamma$ qui fixe $\xi_1$ et \'echange $\xi_{0}$ et $\xi_2$.

Notons $v_1=s_{\frac{1}{2}}s_0v_0=s_{\frac{1}{2}}v_0$. Comme
$s_{\frac{1}{2}}$ envoie la g\'eod\'esique entre $\xi_1$ et $\xi_0$
sur la g\'eod\'esique entre $\xi_1$ et $\xi_2$, l'\'el\'ement
$s_{\frac{1}{2}}s_0$ de $\Gamma$ envoie $x_0$ sur $x_1$
  et $v_0$ sur $v_1$, qui est une ar\^ete issue de $x_1$ qui ne
rentre ni dans $H\!B_{\xi_1}$ ni dans $H\!B_{\xi_2}$.

De m\^eme, notons $s_{-{1\over 2}}=s_{-{1\over 2}}(\ell)$ l'unique
\'el\'ement (d'ordre 2) de $\Gamma$ qui fixe $\xi_0$ et \'echange
$\xi_1$ et $\xi_{-1}$. Posons $v_{-1}=s_{-{1\over 2}} s_0v_0$, c'est
une ar\^ete issue de $x_{-1}$ qui ne rentre ni dans $H\!B_{\xi_0}$, ni
dans $H\!B_{\xi_{-1}}$.

Posons $\wt T(\ell)=(\xi_-,\xi_+,x_1,v_1)$. L'application $\wt
T:\G_0(\TT_q)\ra\G_0(\TT_q)$ est un ho\-m\'e\-o\-mor\-phis\-me (dont
l'inverse est $ \ell\mapsto(\xi_-,\xi_+,x_{-1},v_{-1})$), appel\'e
{\it application de premier retour du flot g\'eod\'esique sur
    $\TT_q$}.  Clairement, l'application $\wt T$ anti-commute avec
l'application de renversement du temps $\wt\tau$~:
$$\wt\tau\circ\wt T\circ \wt\tau= \wt T^{-1}\;.$$
Pour tout $\gamma$ dans $\Gamma$, on a $s_0(\gamma\ell)=\gamma
s_0(\ell)\gamma^{-1}$ et $s_{\frac{1}{2}}(\gamma\ell)=\gamma
s_{\frac{1}{2}}(\ell)\gamma^{-1}$.  Donc $\wt T$ est \'equivariante,
et induit par passage au quotient un hom\'eomorphisme
$$T:\G_0(\Gamma\bac\TT_q)\ra \G_0(\Gamma\bac\TT_q),$$
qui anti-commute avec l'application de renversement du temps $\tau$,
et qui est appel\'e {\it application de premier retour du flot
    g\'eod\'esique sur $\Gamma\bac\TT_q$}.

Notons $\Gamma_\infty'$ l'ensemble des \'el\'ements de
$\Gamma_\infty$ qui ne fixent pas $x_*$. Donc
$\Gamma_\infty'=\Gamma_\infty-({\rm PGL}(2,k)\cap\Gamma_\infty)$.  On
munit $\Gamma'_\infty$ de la topologie discr\`ete et
$(\Gamma'_\infty)^\ZZ$ de la topologie produit.

\medskip
Le r\'esultat suivant est un premier th\'eor\`eme de codage
du flot g\'eod\'esique sur le quotient de l'arbre de Bruhat-Tits de
${\rm PGL}(2,k)$ par le sous-groupe arithm\'etique ${\rm PGL}(2,A)$.
M\^eme si ce n'est pas la version utilis\'ee pour les applications
arithm\'etiques, nous le donnons sous cette forme, qui doit \^etre
g\'en\'eralisable en rempla\c{c}ant $\underline{{\rm PGL}_2}$ par
n'importe quel groupe alg\'ebrique $\underline{G}$ connexe semi-simple
sur $\wh K$, de $\wh K$-rang $1$, et ${\rm PGL}(2,A)$ par n'importe
quel sous-groupe arithm\'etique $\Ga'$ de $\underline{G}(\wh K)$, en
utilisant le th\'eor\`eme de structure du graphe de groupes $\Ga'\bac
T'$, quotient de $T'$ par $\Ga'$, avec $T'$ l'arbre de Bruhat-Tits de
$(\underline{G},\wh K)$.  En effet, le graphe de groupes $\Ga'\bac T'$
est obtenu \`a partir d'un graphe fini de groupes finis, en rajoutant
un nombre fini de rayons de groupes finis, analogues au rayon
modulaire $\Ga\bac\TT_q$ (voir par exemple \cite{Ser2}).

\btheo \label{theo:premretcommutdecal}
Il existe un hom\'eomorphisme $\Theta:\G_0(\Gamma\bac\TT_q)\ra
(\Gamma'_\infty)^\ZZ$ tel que les diagrammes suivants commutent~:
$$\begin{array}{ccc}
\G_0(\Gamma\bac\TT_q) & \stackrel{\Theta}{\longrightarrow} &
(\Gamma'_\infty)^\ZZ \\
\;\downarrow T & & \;\downarrow \sigma \\
\G_0(\Gamma\bac\TT_q) & \stackrel{\Theta}{\longrightarrow} &
(\Gamma'_\infty)^\ZZ
\end{array}
\;\;\;\;\;\;\;
\begin{array}{ccc}
\G_0(\Gamma\bac\TT_q) & \stackrel{\Theta}{\longrightarrow} &
\!\!\!(\Gamma'_\infty)^\ZZ\;\;\;\; \\
\;\downarrow \tau & & \;\;\;\;\downarrow \kappa\circ\sigma \\
\G_0(\Gamma\bac\TT_q) & \stackrel{\Theta}{\longrightarrow} &
\!\!\!(\Gamma'_\infty)^\ZZ \;\;\;\;
\end{array}$$
avec $\sigma$ le d\'ecalage \`a gauche et $\kappa$
la transformation  de $(\Gamma'_\infty)^\ZZ$ d\'efinie par
$\kappa((\beta_n)_{n\in\ZZ})= (\beta'_n)_{n\in\ZZ}$ avec 
$\beta'_{n}=\beta_{-n}^{-1}.$
\etheo

\begin{center}
\begin{picture}(0,0)%
\includegraphics{fig_notationcodage.pstex}%
\end{picture}%
\setlength{\unitlength}{3947sp}%
\begingroup\makeatletter\ifx\SetFigFont\undefined%
\gdef\SetFigFont#1#2#3#4#5{%
  \reset@font\fontsize{#1}{#2pt}%
  \fontfamily{#3}\fontseries{#4}\fontshape{#5}%
  \selectfont}%
\fi\endgroup%
\begin{picture}(6848,2895)(748,-3451)
\put(7134,-1123){\makebox(0,0)[lb]{\smash{\SetFigFont{12}{14.4}{\rmdefault}{\mddefault}{\updefault}{$0$}%
}}}
\put(3141,-2274){\makebox(0,0)[lb]{\smash{\SetFigFont{12}{14.4}{\rmdefault}{\mddefault}{\updefault}{$x_1$}%
}}}
\put(1266,-2941){\makebox(0,0)[lb]{\smash{\SetFigFont{12}{14.4}{\rmdefault}{\mddefault}{\updefault}{$x_{-2}$}%
}}}
\put(1791,-2266){\makebox(0,0)[lb]{\smash{\SetFigFont{12}{14.4}{\rmdefault}{\mddefault}{\updefault}{$x_{-1}$}%
}}}
\put(4059,-1421){\makebox(0,0)[lb]{\smash{\SetFigFont{12}{14.4}{\rmdefault}{\mddefault}{\updefault}{$\eta_1$}%
}}}
\put(2429,-1966){\makebox(0,0)[lb]{\smash{\SetFigFont{12}{14.4}{\rmdefault}{\mddefault}{\updefault}{$x_0$}%
}}}
\put(4416,-2211){\makebox(0,0)[lb]{\smash{\SetFigFont{12}{14.4}{\rmdefault}{\mddefault}{\updefault}{$\xi_2=\xi_{-1}(\wt\tau\ell)$}%
}}}
\put(6459,-1933){\makebox(0,0)[lb]{\smash{\SetFigFont{12}{14.4}{\rmdefault}{\mddefault}{\updefault}{$x_*$}%
}}}
\put(1369,-926){\makebox(0,0)[lb]{\smash{\SetFigFont{12}{14.4}{\rmdefault}{\mddefault}{\updefault}{$\xi_0=\xi_1(\wt\tau\ell)$}%
}}}
\put(6489,-973){\makebox(0,0)[lb]{\smash{\SetFigFont{12}{14.4}{\rmdefault}{\mddefault}{\updefault}{$1$}%
}}}
\put(5709,-1168){\makebox(0,0)[lb]{\smash{\SetFigFont{12}{14.4}{\rmdefault}{\mddefault}{\updefault}{$\infty$}%
}}}
\put(3636,-2941){\makebox(0,0)[lb]{\smash{\SetFigFont{12}{14.4}{\rmdefault}{\mddefault}{\updefault}{$x_2$}%
}}}
\put(3338,-956){\makebox(0,0)[lb]{\smash{\SetFigFont{12}{14.4}{\rmdefault}{\mddefault}{\updefault}{$\xi_1$}%
}}}
\put(7200,-2821){\makebox(0,0)[lb]{\smash{\SetFigFont{12}{14.4}{\rmdefault}{\mddefault}{\updefault}{$\ga_0^{-1}\xi_+$}%
}}}
\put(2379,-736){\makebox(0,0)[lb]{\smash{\SetFigFont{12}{14.4}{\rmdefault}{\mddefault}{\updefault}{$\eta_0=\eta_0(\wt\tau \ell)$}%
}}}
\put(751,-1298){\makebox(0,0)[lb]{\smash{\SetFigFont{12}{14.4}{\rmdefault}{\mddefault}{\updefault}{$\eta_{-1}$}%
}}}
\put(846,-3451){\makebox(0,0)[lb]{\smash{\SetFigFont{12}{14.4}{\rmdefault}{\mddefault}{\updefault}{$\xi_-$}%
}}}
\put(4041,-3429){\makebox(0,0)[lb]{\smash{\SetFigFont{12}{14.4}{\rmdefault}{\mddefault}{\updefault}{$\xi_+$}%
}}}
\put(5363,-2851){\makebox(0,0)[lb]{\smash{\SetFigFont{12}{14.4}{\rmdefault}{\mddefault}{\updefault}{$\ga_0^{-1}\xi_-$}%
}}}
\put(4989,-1648){\makebox(0,0)[lb]{\smash{\SetFigFont{12}{14.4}{\rmdefault}{\mddefault}{\updefault}{$\ga_0$}%
}}}
\put(3641,-1826){\makebox(0,0)[lb]{\smash{\SetFigFont{12}{14.4}{\rmdefault}{\mddefault}{\updefault}{$v_1$}%
}}}
\put(6481,-1423){\makebox(0,0)[lb]{\smash{\SetFigFont{12}{14.4}{\rmdefault}{\mddefault}{\updefault}{$v_*$}%
}}}
\put(2479,-1311){\makebox(0,0)[lb]{\smash{\SetFigFont{12}{14.4}{\rmdefault}{\mddefault}{\updefault}{$v_0$}%
}}}
\put(1136,-1846){\makebox(0,0)[lb]{\smash{\SetFigFont{12}{14.4}{\rmdefault}{\mddefault}{\updefault}{$v_{-1}$}%
}}}
\put(4316,-2826){\makebox(0,0)[lb]{\smash{\SetFigFont{12}{14.4}{\rmdefault}{\mddefault}{\updefault}{$v_2$}%
}}}
\end{picture}
\end{center}
\begin{center}{\bf Figure \addtocounter{fig}{1}\arabic{fig} :}
Codage des g\'eod\'esiques.
\end{center}

\dem
Pour tout $n$ dans $\ZZ$ et pour toute g\'eod\'esique $\ell$ de
$\G_0(\TT_q)$, notons $\eta_n=\eta_n(\ell)=\eta_0(\wt T^n\ell)$ et
$v_n=v_n(\ell)=v_0(\wt T^n\ell)$. Par d\'efinition de $\wt T$, les
notations $v_1$ et $v_{-1}$ co\"{\i}ncident avec celles d\'ej\`a
introduites.  Soit $\gamma_n=\gamma_n(\ell)$ l'unique \'el\'ement de
$\Gamma$ envoyant $\infty,1,0$ sur respectivement
$\xi_n,\eta_n,\xi_{n+1}$.  Notons $i$ l'unique \'el\'ement (d'ordre
$2$) de $\Gamma$ fixant $1$ et \'echangeant $0$ et $\infty$. Posons
$$\beta_n=\beta_n(\ell)=i\gamma_{n-1}^{-1}\gamma_{n}\;\;\;\;{\rm
     et}\;\;\;\;\wt\Theta(\ell)=(\beta_n)_{n\in\ZZ}\;.
$$

Comme $i\gamma_{n-1}^{-1}\gamma_{n}\infty=i\gamma_{n-1}^{-1}\xi_{n}=i0
=\infty$, l'\'el\'ement $\beta_n$ de $\Gamma$ appartient \`a
$\Gamma_\infty$.  De plus, $\beta_n(x_*)= i\gamma_{n-1}^{-1}
\gamma_{n}(x_*) = i\gamma_{n-1}^{-1}(x_n )$, donc $\beta_n(x_*)$ ne
peut pas valoir $x_*$, car $x_{n-1}$ est diff\'erent de $x_n$ et
$i\gamma_{n-1}^{-1}(x_{n-1})=x_*$.  Ainsi $\beta_n$ est dans
$\Gamma_\infty'$.

Pour tout $\gamma$ de $\Gamma$, on a la relation de naturalit\'e
$\gamma_n(\gamma\ell)=\gamma\gamma_n(\ell)$.  L'application
$\wt\Theta:\G_0(\TT_q)\ra (\Gamma'_\infty)^\ZZ$ est constante sur
chaque orbite de $\Gamma$ dans $\G_0(\TT_q)$. Elle induit donc par
passage au quotient une application $\Theta:
\G_0(\Gamma\bac\TT_q)\ra(\Gamma'_\infty)^\ZZ$.

\medskip
\noindent{\bf Commutativit\'e des diagrammes.}  Par construction de
$\wt T$ et par d\'efinition des $\eta_n$, pour tout $n$ dans $\ZZ$, on
a $\xi_n(\ell)=\xi_0(\wt T^n\ell)$ et $\eta_n(\ell)=\eta_0(\wt
T^n\ell)$.  Donc $$\xi_n(\wt T\ell)=\xi_{n+1}(\ell)\;\;\;{\rm
    et}\;\;\;\eta_n(\wt T\ell)=\eta_{n+1}(\ell)\;.$$
Par cons\'equent, $\ga_n(\wt T\ell)=\ga_{n+1}(\ell)$, d'o\`u
$\beta_n(\wt T\ell)=\beta_{n+1}(\ell)$.  Donc $\Theta\circ
T=\sigma\circ\Theta$, ce qui montre la commutativit\'e du diagramme de
gauche du th\'eor\`eme \ref{theo:premretcommutdecal}.

Par d\'efinition des $\xi_n$, on a $\xi_n(\wt\tau
\ell)=\xi_{-(n-1)}(\ell)$. Comme $\eta_0(\wt\tau\ell)=\eta_0(\ell)$,
et comme $\wt T$ et $\wt \tau$ anti-commutent, par d\'efinition des
$\eta_n$, on a $\eta_n(\wt\tau\ell)=\eta_{-n}(\ell)$. La d\'efinition
des $\gamma_n$ montre que $\gamma_n(\wt\tau\ell)$ et
$\gamma_{-n}(\ell)i$ co\"{\i}ncident en $\infty,1,0$. Donc
$\gamma_n(\wt\tau\ell)=\gamma_{-n}(\ell)i$. D'o\`u
$$\beta_n(\wt\tau\ell)=i\gamma_{n-1}^{-1}(\wt\tau\ell)
\gamma_{n}(\wt\tau\ell) =\ga_{-(n-1)}^{-1}(\ell)\ga_{-n}(\ell)i=
\beta_{-(n-1)}^{-1}(\ell)\;.$$
Ceci montre la commutativit\'e du diagramme de droite du th\'eor\`eme
\ref{theo:premretcommutdecal}.

\medskip
\noindent{\bf Continuit\'e de $\Theta$.}  Montrons que l'application
$\ell\mapsto\beta_1(\ell)$ est localement constante.  Soit $\ell$ un
\'el\'ement de $\G_0(\TT_q)$. Si $\ell'$ est une g\'eod\'esique
d\'ecor\'ee, ayant m\^emes origine et d\'ecoration que $\ell$, et
co\"{\i}ncidant avec $\ell$ sur le $1$-voisinage du segment
$[x_0(\ell),x_1(\ell)]$, alors $\eta_0(\ell)=\eta_0(\ell')$ et
$\xi_k(\ell)=\xi_k(\ell')$ pour $k=0,1,2$. Donc
$\eta_1(\ell)=\eta_1(\ell')$ par d\'efinition de $\eta_1$.  Par
d\'efinition des $\ga_n$, on a alors $\ga_k(\ell)=\ga_k(\ell')$ pour
$k=0,1$. Par cons\'equent $\beta_1(\ell)=\beta_1(\ell')$.  Donc
$\ell\mapsto\beta_1(\ell)$ est continue. Comme
$\beta_n(\ell)=\beta_1(\wt T^{n-1}\ell)$ et puisque $\wt T$ est un
hom\'eomorphisme, l'application $\ell\mapsto\beta_n(\ell)$ est
continue sur $\G_0(\TT_q)$ pour tout $n$ dans $\ZZ$. Donc $\wt\Theta$
est continue sur $\G_0(\TT_q)$, et par passage au quotient, $\Theta$
est continue sur $\G_0(\Gamma\bac\TT_q)$.

\medskip
\noindent{\bf Injectivit\'e de $\Theta$.}  Soient $\ell$ et $\ell'$ deux
\'el\'ementss de $\G_0(\TT_q)$ tels que $\wt\Theta(\ell)=
\wt\Theta(\ell')$.  Montrons que $\ell $ et $\ell'$ sont dans la
m\^eme orbite sous $\Gamma$. Ceci montrera l'injectivit\'e de
$\Theta$. Quitte \`a remplacer $\ell$ par $\ga_0(\ell)^{-1}\ell$ et
$\ell'$ par $\ga_0(\ell')^{-1}\ell'$, on peut supposer que
$\xi_0(\ell)=\xi_0(\ell')=\infty$, $\eta_0(\ell)=\eta_0(\ell')=1$ et
$\xi_1(\ell)=\xi_1(\ell')=0$. En particulier, $\ell$ et $\ell'$ ont
m\^emes origine et d\'ecoration. De plus, $\ga_0(\ell)=\ga_0(\ell')=
{\rm id}$.  Pour $n\geq 0$, on a
$$\gamma_n=\gamma_0(\gamma_0^{-1}\gamma_1)\cdots
(\gamma_{n-1}^{-1}\gamma_n)= \ga_0i\beta_1i\beta_2\cdots i\beta_n$$
et
$$\ga_{-n}=\ga_0(\ga_{-1}^{-1}\ga_0)^{-1}(\ga_{-2}^{-1}\ga_{-1})^{-1}
\cdots(\ga_{-n}^{-1}\ga_{-(n-1)})^{-1}=
\ga_0\beta_{0}^{-1}i\beta_{-1}^{-1}i\cdots \beta_{-(n-1)}^{-1}i\;.$$
Donc pour tout $n$ dans $\ZZ$, on a $\gamma_n(\ell)=\gamma_n(\ell')$.
Comme $\ga_n\infty=\xi_n$, les points $\xi_n(\ell)$ et $\xi_n(\ell')$
co{\"\i}ncident. Or deux g\'eod\'esiques d\'ecor\'ees ayant m\^emes
origine et d\'ecoration et traversant les m\^emes horoboules de $\H\B$
sont \'egales. Donc $\ell=\ell'$.

\medskip
\noindent{\bf Surjectivit\'e de $\Theta$.}  Soit $(\beta_n)_{n\in\ZZ}$
dans $(\Gamma'_\infty)^\ZZ$. Posons $\ga_0={\rm id}$. Par r\'ecurrence
pour tout $n$ dans $\ZZ$, posons $\gamma_n=\gamma_{n-1}i\beta_n$.
Notons $x_n=\gamma_nx_\ast$.  Comme les $\beta_n$ sont dans
$\Gamma_\infty'$, les points $x_n$ et $x_{n+1}$ sont diff\'erents.
Consid\'erons la courbe $\ell$, g\'eod\'esique par morceaux, obtenue
en recollant cons\'ecutivement les segments $[x_n,x_{n+1}]$ (qui ne
sont pas r\'eduits \`a des points). Posons
$$\xi_n=\ga_n\infty\;\;\;{\rm et}\;\;\;\eta_n=\ga_n 1\;,$$
qui sont
des points de $\PP_1(K)$.  Alors
$\xi_{n+1}=\ga_{n+1}\infty=\ga_{n}i\beta_n\infty=\ga_n0.$ Comme
$\{x_*\}=H\!B_\infty\cap H\!B_1\cap H\!B_0$, on a
$\{x_n\}=H\!B_{\xi_n}\cap H\!B_{\eta_n}\cap H\!B_{\xi_{n+1}}$.  Par
convexit\'e, le segment $[x_n,x_{n+1}]$ est contenu dans l'horoboule
$H\!B_{\xi_{n+1}}$. Comme deux horoboules distinctes de $\H\B$ sont
d'int\'erieurs disjoints, il vient donc $[x_n,x_{n+1}]\cap
[x_{n+1},x_{n+2}]=\{x_{n+1}\}$. Par cons\'equent $\ell$ est une
g\'eod\'esique. Il est alors imm\'ediat que, munie de la d\'ecoration
$v_\ast$, la g\'eod\'esique $\ell$ a pour image $(\beta_n)_{n\in\ZZ}$
par $\wt\Theta$. Donc $\wt\Theta$ est surjective, et, par passage au
quotient, $\Theta$ aussi. \cqfd

\subsection{Lien avec la transformation d'Artin}
\label{subsec:lienartin}

Notons $\G'_0(\TT_q)$ l'ensemble des g\'eod\'esiques d\'ecor\'ees
de $\G_0(\TT_q)$ qui ont pour origine $x_*$ et pour d\'ecoration
$v_*$.

\blemm
Soit $\ell$ un \'el\'ement de $\G_0(\TT_q)$, alors $\ell$ appartient
\`a $\G'_0(\TT_q)$ si et seulement si $\ga_0(\ell)={\rm id}$.
\elemm

\dem
Si $\ell$ est dans $\G'_0(\TT_q)$, alors par d\'efinition de
$\ga_0$, on a $\ga_0(\ell)={\rm id}$. R\'e\-ci\-pro\-que\-ment, si
$\ga_0(\ell)={\rm id}$, alors $\xi_0(\ell)=\infty$, $\xi_1(\ell)=0$ et
$\eta_0(\ell)=1$. Comme $x_0(\ell)$ est le point d'intersection des
horosph\`eres $H_\infty$ et $H_0$, il n'est autre que $x_*$. De
m\^eme, $v_0(\ell)$ est l'ar\^ete issue de $x_*$ qui rentre dans
$H\!B_1$, donc c'est $v_*$.  Donc $\ell$ appartient \`a
$\G'_0(\TT_q)$.
\cqfd

\blemm
L'application $\varphi$ de $\G_0(\Gamma\bac\TT_q)$ dans
$\G_0'(\TT_q)$, d\'efinie par $\varphi(\overline
\ell)=\gamma_0(\ell)^{-1}\ell$ o\`u $\ell$ est un repr\'esentant de
l'\'el\'ement $\overline \ell$ de $\G_0(\Gamma\bac\TT_q)$, est un
hom\'eomorphisme, qui est une section de la projection canonique
$\pi_0:\G_0(\TT_q)\ra \G_0(\Gamma\bac\TT_q)$, que nous appelerons la
{\rm section canonique}.
\elemm

\dem
Si $\ell$ est une g\'eod\'esique d\'ecor\'ee de $\G_0(\TT_q)$,
alors $\ell'=\ga_0(\ell)^{-1}\ell$ est dans $\G'_0(\TT_q)$, car
$\ga_0(\ell')={\rm id}$.  Comme deux g\'eod\'esiques d\'ecor\'ees qui
ont m\^eme d\'ecoration et qui co\"{\i}ncident entre $-1$ et $1$ ont
m\^eme $\ga_0$, l'application $\tilde{\varphi}$ de $\G_0(\TT_q)$ dans
$\G_0'(\TT_q)$ d\'efinie par $\ell\mapsto\ga_0(\ell)^{-1}\ell$ est
continue.  Elle est surjective car $\ga_0(\ell)={\rm id}$ si $\ell$
est dans $\G'_0(\TT_q)$. Comme $\ga_0(\ga \ell)=\ga \ga_0(\ell)$ pour
tout $\ga$ dans $\Ga$, deux \'el\'ements $\ell,\ell'$ de $\G_0(\TT_q)$
sont dans la m\^eme $\Gamma$-orbite si et seulement si
$\ga_0(\ell)^{-1}\ell=\ga_0(\ell')^{-1}\ell'$. L'application
$\tilde{\varphi}$ induit donc par passage au quotient une application
continue bijective $\varphi:\G_0(\Gamma\bac\TT_q) \ra \G_0'(\TT_q)$.
Comme $\ga_0(\ell)^{-1} \ell$ est dans la m\^eme $\Ga$-orbite que
$\ell$, on a $\pi_0\circ\varphi={\rm id}$. Si
$j:\G'_0(\TT_q)\ra\G_0(\TT_q)$ est l'inclusion, alors
$\varphi^{-1}=\pi_0\circ j$, donc $\varphi^{-1}$ est continue.
\cqfd

\medskip
Posons $$J=\bigcup_{a\in A-k}(a+X^{-1}\O)\;.$$
Si $\ell$ est dans $\G_0'(\TT_q)$, alors ses extr\'emit\'es
$\xi_-=\xi_-(\ell)$ et $\xi_+=\xi_+(\ell)$ appartiennent
respectivement \`a $J\cap \,^c\!K$ et $X^{-1}\O\cap \,^c\!K$ (voir par
exemple \cite{Pau} et la figure 1).  L'application de $\G_0'(\TT_q)$
dans $(J\cap \,^c\!K)\times (X^{-1}\O\cap \,^c\!K)$ d\'efinie par
$\ell\mapsto (\xi_-,\xi_+)$ est un hom\'eomorphisme. Nous
identifierons dans la suite une g\'eod\'esique $\ell$ de $
\G'_0(\TT_q)$ avec son couple d'extr\'emit\'es irrationnelles
$(\xi_-,\xi_+)$ dans $ (J\cap \,^c\!K)\times (X^{-1}\O\cap \,^c\!K)$.

\medskip Nous allons modifier le codage obtenu dans le th\'eor\`eme
\ref{theo:premretcommutdecal} pour obtenir un nouveau codage, reli\'e
cette fois-ci \`a la transformation d'Artin $\Psi$. Pour cela nous
allons pr\'eciser les homographies $\gamma_n$ qui sont apparues au
cours de la d\'emonstration pr\'ec\'edente. On introduit les notations
suivantes~: on appelle $i$ (inversion de centre 0), $t_a$ (translation
de $a$) et $\lambda_\alpha$ (homoth\'etie de rapport $\alpha^2$) les
homographies associ\'ees aux matrices $\left[\begin{array}{cc}
      0&1\\1&0\end{array}\right]$, $\left[\begin{array}{cc}
      1&a\\0&1\end{array}\right]$ et $\left[\begin{array}{cc}
      \alpha&0\\0&\alpha^{-1}\end{array}\right]$ avec $a$ dans $A$ et
$\alpha$ dans $k^\times$.

Tout \'el\'ement $\beta$ de $\Gamma_\infty$ s'\'ecrit alors de
mani\`ere unique sous la forme $\beta=t_a\lambda_\alpha$ avec $a$ dans
$A$ et $\alpha$ dans $k^\times$. Si $\beta$ est de plus dans
$\Gamma'_\infty$ alors $a$ ne peut pas \^etre dans $k$ sinon $\beta$
serait un \'el\'ement de ${\rm PGL}(2,k)$. Donc tout \'el\'ement
$\beta$ de $\Gamma_\infty'$ peut s'\'ecrire de fa\c con unique sous la
forme $\beta=t_a\lambda_\alpha$ avec $a$ dans $A-k$ et $\alpha$ dans
$k^\times$.

On remarque qu'on a les relations \'el\'ementaires suivantes~:
$t_a^{-1}=t_{-a}$, $\lambda_\alpha^{-1}=\lambda_{1/\alpha}$,
$t_a\lambda_\alpha=\lambda_{\alpha}t_{a/\alpha^2}$
et $i\lambda_\alpha=\lambda_{1/\alpha}i$.

Montrons maintenant  le r\'esultat suivant~:

\begin{prop}\label{prop:techniqueavantcodagedeux}
     Soit $\ell=(\xi_-,\xi_+)$ une g\'eod\'esique de $\G'_0(\TT_q)$. On
     note $(a_n)_{n\geq 1}$ le d\'e\-ve\-lop\-pe\-ment en fractions
continues d'Artin de
     $\xi_+$ et $(a_{-n})_{n\geq 0}$ celui de ${-1\over \xi_-}$. Alors pour
     tout entier $n>0$~:

     il existe $\alpha_n$ dans $k^\times$ tel que
     $\gamma_n(\ell)=it_{a_1}it_{a_2}\dots it_{a_n}\lambda_{\alpha_n}$,

     il existe $\alpha_{-n}$ dans $k^\times$ tel que
     $\gamma_{-n}(\ell)=t_{-a_0}it_{-a_{-1}}\dots
     it_{-a_{-(n-1)}}i\lambda_{\alpha_{-(n-1)}}$.
\end{prop}

\dem
Soit $\ell=(\xi_-,\xi_+)$ une g\'eod\'esique de $\G'_0(\TT_q)$.
Proc\`edons par r\'ecurrence sur $n$. Comme $\ga_0= {\rm id}$, la
d\'emonstration du th\'eor\`eme \ref{theo:premretcommutdecal} montre
qu'il existe une suite $(\beta_n)_{n\in\ZZ}$ d'\'el\'ements de
$\Gamma'_\infty$ telle que pour tout $n>0$, on a~:
$$\gamma_n=i\beta_1i\beta_2\dots i\beta_n \;\mbox{ et }\;
\gamma_{-n}=\beta_0^{-1}i\beta_{-1}^{-1}i\dots
\beta_{-(n-1)}^{-1}i\;.$$
La transformation $\beta_1$ est dans $\Gamma'_\infty$. Elle s'\'ecrit
donc de mani\`ere unique sous la forme $t_{a_1}\lambda_{\alpha_1}$
avec $a_1$ dans $A-k$ et $\alpha_1$ dans $k^\times$. Nous avons donc
bien $\gamma_1=it_{a_1}\lambda_{\alpha_1}$.  Supposons que $\gamma_n$
s'\'ecrive $it_{a_1}it_{a_2}\dots it_{a_n}\lambda_{\alpha_n}$. Alors
$$
\gamma_{n+1}=\gamma_n i\beta_{n+1}=it_{a_1}it_{a_2}\dots
it_{a_n}\lambda_{\alpha_n}i\beta_{n+1}=it_{a_1}it_{a_2}\dots
it_{a_n}i\lambda_{1/\alpha_n}\beta_{n+1}\;.
$$
La transformation $\lambda_{1/\alpha_n}\beta_{n+1}$ est dans
$\Gamma'_\infty$, elle s'\'ecrit donc sous la forme
$t_{a_{n+1}}\lambda_{\alpha_{n+1}}$ avec $a_{n+1}$ dans $A-k$ et
$\alpha_{n+1}$ dans $k^\times$. Ainsi $\gamma_{n+1}=
it_{a_1}it_{a_2}\dots it_{a_n}i t_{a_{n+1}}\lambda_{\alpha_{n+1}}$.

De fa\c con analogue, traitons les $\gamma_{-n}$. On a
$\gamma_{-1}=\beta_0^{-1}i$. La transformation $\beta_0$ est dans
$\Gamma'_\infty$, elle s'\'ecrit donc de mani\`ere unique sous la
forme $\lambda_{\alpha_0}t_{a_0}$ et donc $\gamma_{-1}=
(\lambda_{\alpha_0}t_{a_0})^{-1}i=t_{-a_0}
\lambda_{1/\alpha_0}i=t_{-a_0}i\lambda_{\alpha_0}$. Supposons que
$$\gamma_{-n}=t_{-a_0}it_{-a_{-1}}\dots it_{-a_{-(n-1)}}
i\lambda_{\alpha_{-(n-1)}}\;.$$ Alors,
$$
\gamma_{-(n+1)}=\gamma_{-n}\beta_{-n}^{-1}i=t_{-a_0}it_{-a_{-1}}\dots
it_{-a_{-(n-1)}}i\lambda_{\alpha_{-(n-1)}}\beta_{-n}^{-1}i\;.$$
La transformation $(\lambda_{\alpha_{-(n-1)}}\beta_{-n}^{-1})^{-1}=
\beta_{-n}\lambda_{1/\alpha _{-(n-1)}}$ est dans
$\Gamma'_\infty$. Elle s'\'ecrit donc de mani\`ere unique sous la
forme $\lambda_{\alpha_{-n}}t_{a_{-n}}$. Alors
\begin{eqnarray*}
\gamma_{-(n+1)} & =&t_{-a_0}it_{-a_{-1}}\dots
it_{-a_{-(n-1)}}i(\lambda_{\alpha_{-n}}t_{a_{-n}})^{-1}i\\
&=& t_{-a_0}it_{-a_{-1}}\dots
it_{-a_{-(n-1)}}it_{-a_{-n}}\lambda_{1/\alpha_{-n}}i\\
&=&
t_{-a_0}it_{-a_{-1}}\dots
it_{-a_{-(n-1)}}it_{-a_{-n}}i\lambda_{\alpha_{-n}}.
\end{eqnarray*}
Par cons\'equent, pour tout entier $n>0$,
$$
\xi_n=\gamma_n\infty=it_{a_1}it_{a_2}\dots
it_{a_n}\lambda_{\alpha_n}\infty=it_{a_1}it_{a_2}\dots it_{a_{n-1}} 0
$$
et
$$
\xi_{-n}=\gamma_{-n}\infty=t_{-a_0}it_{-a_{-1}}\dots
it_{-a_{-(n-1)}}i\lambda_{\alpha_{-(n-1)}}\infty=-t_{a_0}
it_{a_{-1}}\dots it_{a_{-(n-1)}}0\;.
$$
Comme $\xi_+$ est la limite quand $n$ tend vers l'infini de la suite
$\xi_n$, le d\'eveloppement en fractions continues d'Artin de
$\xi_+$ est $(a_n)_{n\geq 1}$.

De m\^eme $\xi_-$ est la limite de la suite $\xi_{-n}$, donc ${-1\over
    \xi_- }=-i\xi_-$ est la limite de la suite $-i\xi_n=
it_{a_0}it_{a_{-1}}\dots it_{a_{-(n-1)}}0$. Le d\'eveloppement en
fractions continues d'Artin de ${-1\over \xi_- }$ est donc
$(a_{-n})_{n\geq 0}$.  \cqfd

\bigskip D\'efinissons la transformation
$\wt\Psi:\G_0'(\TT_q)\ra\G_0'(\TT_q)$ par
$\wt\Psi(\ell)=\gamma_1^{-1}(\ell)\wt T(\ell)$. Il est imm\'ediat
que le diagramme suivant commute~:
$$\begin{array}[b]{ccc}
\G_0(\Gamma\bac\TT_q) & \stackrel{T}{\longrightarrow} &
\G_0(\Gamma\bac\TT_q)
\\
\downarrow \varphi & & \downarrow \varphi \\
\G_0'(\TT_q) &\stackrel{\wt\Psi}{\longrightarrow} &
\G_0'(\TT_q)
\end{array}.$$
Ceci revient \`a dire que $\wt\Psi$ est conjugu\'ee par la section
canonique $\varphi$ \`a l'application $T$ de premier retour du flot
g\'eod\'esique sur $\G_0(\Gamma\bac\TT_q)$.

\medskip
Montrons maintenant le r\'esultat suivant, qui implique le
th\'eor\`eme \ref{theo:intro} de l'introduction~:

\btheo\label{theo:codagedeux}
Soit $\ell=(\xi_-,\xi_+)$ un \'el\'ement de $\G_0'(\TT_q)$, soit
$(a_n)_{n\geq 1}$ le d\'eveloppement en fractions continues d'Artin de
$\xi_+$ et $(a_{-n})_{n\geq 0}$ celui de ${-1\over \xi_-}$. Notons
$\Theta'$ l'application de $\G_0'(\TT_q)$ dans $(A-k)^{\ZZ}$ d\'efinie
par $\Theta'(\ell)=(a_n)_{n\in\ZZ}$.  Alors $\Theta'$ est un
hom\'eomorphisme qui rend le diagramme suivant commutatif~:
$$\begin{array}[b]{ccc} \G_0'(\TT_q) &
    \stackrel{\Theta'}{\longrightarrow} & (A-k)^{\ZZ}
    \\
    \wt\Psi\downarrow\;  & & \downarrow \sigma \\
    \G_0'(\TT_q) &\stackrel{\Theta'}{\longrightarrow} & (A-k)^{\ZZ}
\end{array}\;,$$
o\`u $\sigma$ est le d\'ecalage \`a gauche des suites bilat\`eres de
$(A-k)^{\ZZ}$. De plus~:
$$\wt\Psi(\xi_-,\xi_+)=\left({1\over \xi_-}-\left[{1\over
       \xi_+}\right],\Psi(\xi_+)\right)\;,$$
o\`u $\Psi$ est la transformation d'Artin.
\etheo

\dem Comme un \'el\'ement de $X^{-1}\O\cap\,^c\!K$ est uniquement
d\'etermin\'e par son d\'eveloppement (infini) en fractions continues,
et comme $(\xi_-,\xi_+)$ est dans
$(J\cap\,^c\!K)\times(X^{-1}\O\cap\,^c\!K)$, alors ${-1\over \xi_-}$
est dans $X^{-1}\O\cap\,^c\!K$ et donc l'application $\Theta'$ est une
bijection. Comme l'application d'Artin est continue (car localement
constante) sur $X^{-1}\O\cap\,^c\!K$, et par d\'efinition des $a_n$,
l'application $\Theta'$ est continue.  D'apr\`es la d\'emonstration
de la proposition \ref{prop:techniqueavantcodagedeux}, si $\ell$ et
$\ell'$ sont deux \'el\'ements de $\G_0'(\TT_q)$ qui ont m\^emes $a_n$
pour $-N\leq n\leq N$, alors $\ell$ et $\ell'$ ont m\^emes $\xi_n$
pour $-(N-1)\leq n\leq N-1$, donc co\"{\i}ncident au moins sur
$-(N-1)\leq n\leq N-1$.  Donc $\Theta'^{-1}$ est continue.  Par la
proposition pr\'ec\'edente, si $\ell=(\xi_-,\xi_+)$ est dans $
\G_0'(\TT_q)$, alors $\ga_1=it_{a_1}$ o\`u $a_1$ est le premier terme
du d\'eveloppement en fractions continues de $\xi_+$, c'est-\`a-dire
la partie enti\`ere de ${1\over\xi_+}$.  En particulier, le diagramme
de l'\'enonc\'e est \'evidemment commutatif.

Enfin,
$$\wt\Psi(\xi_-,\xi_+)=\ga_1^{-1}(\xi_-,\xi_+)=
((it_{a_1})^{-1}\xi_-,(it_{a_1})^{-1}\xi_+) =
({1\over\xi_-}-a_1,{1\over \xi_+}-a_1)\;.$$
Comme $\Psi(\xi_+)={1\over \xi_+}-a_1$, ceci montre le r\'esultat.
\cqfd

\section{Mesure invariante par le flot g\'eod\'esique sur
le rayon mo\-du\-lai\-re}
\label{sec:mesure}

Notons $\mu=\mu_{\rm Haar}$ la mesure de Haar sur le groupe
topologique additif $\wh{K}$ (identifi\'e avec $\partial
\TT_q-\{\infty\}$), normalis\'ee pour que $\mu(\O)=1$.  Notons $dx_0$
et $dv_0$ les mesures de comptages sur les espaces discrets $V\TT_q$
et $E\TT_q$.

Consid\'erons la mesure $\wt m$ sur l'espace $\G_0(\TT_q)$ d\'efinie,
en utilisant le param\'etrage de $\G_0(\TT_q)$, par
$$d\wt m (\xi_-,\xi_+, x_0,v_0)=
\frac{d\mu(\xi_-)d\mu(\xi_+)dx_0dv_0}{|\xi_+ -\xi_-|^2_\infty}\;.$$

Cette partie est consacr\'e \`a la d\'emonstration du r\'esultat
suivant~:

\btheo \label{theo:mesinva}
La mesure $\wt m$ est une mesure bor\'elienne positive, invariante par
$\Ga$, par $\wt T$ et par renversement du temps $\wt\tau$.  La mesure
$q^2\wt m$ induit par passage au quotient par $\Ga$ une mesure de
probabilit\'e $m$ sur $\G_0(\Ga\bac\TT_q)$, invariante par $T$ et par
renversement du temps $\tau$.
\etheo

Outre une preuve directe, la partie nouvelle de ce théorème est le
calcul de la masse totale de $m$. Mais une fois rappelé les
propositions \ref{prop:etape_un} et \ref{prop:etape_deux}, l'existence
de $m$ découle essentiellement de la construction classique de la
mesure de Bowen-Margulis sur l'espace des géodésiques de $\TT_q$ (voir
par exemple \cite{Coo}), et de \cite{BM} par exemple, qui montre que
les mesures de Patterson-Sullivan et de Hausdorff sur $\partial \TT_q$
co\"{\i}ncident, m\^eme pour les réseaux non uniformes dans
$Aut(\TT_q)$).

\medskip
Dans ce qui suit, l'espace $\wh K$ est muni de la distance $d_\infty$.
Notons $\delta$ la dimension de Hausdorff de $\wh K$, et $\mu_{\rm
     Haus}$ la mesure de Hausdorff (de dimension $\delta$) sur $\wh K$.
Rappelons que pour tout bor\'elien $E$ de $\wh K$, si
$0<\delta<\infty$, nous avons
$$\mu_{\rm Haus}(E)=\lim_{\varepsilon\ra 0}\;\inf \{\sum_{i\in \NN}
r_i^{\delta}\}$$
o\`u la borne inf\'erieure est prise sur tous les recouvrements de $E$
par des boules de rayon $r_i\leq \varepsilon$ pour la distance
$d_\infty$.

Comme $d_\infty$ est invariante par $\Ga_\infty$, il en est de m\^eme
pour $\mu_{\rm Haus}$.

Les deux propositions \ref{prop:etape_un} et \ref{prop:etape_deux}
suivantes sont bien connues, voir par exemple \cite[page 69]{Spr}
pour une d\'emonstration de la seconde.

\bprop\label{prop:etape_un}
La dimension de Hausdorff de $(\wh K,d_\infty)$ est $\delta =1$. La
mesure de Hausdorff de $\O$ est $1$.
\eprop

\dem
Par invariance par translation de $d_\infty$, il suffit de
montrer que la dimension de Hausdorff de $\O$ est $1$, et que sa
mesure de Hausdorff (en dimension $1$) vaut $1$.  L'espace de Cantor
$\O$ est la boule de rayon $1$ et de centre n'importe quel point de
$\O$. Il s'\'ecrit comme l'union disjointe des $q$ parties
$\O_\alpha=\alpha +X^{-1}\O$ pour $\alpha$ dans $k$. Les parties
$\O_\alpha$ sont d'ailleurs les boules de rayon $\frac{1}{q}$ et de
centre n'importe quel point de $\O_\alpha$.  Pour tout $\alpha$ de 
$k$, l'homographie de matrice
$\left[\begin{array}{cc} 1 & \alpha \\ 0 &
      1\end{array}\right]\left[\begin{array}{cc} X^{-1} & 0 \\ 0 &
      1\end{array}\right]$ envoie $\O$ sur $\O_\alpha$, c'est de plus une
  homoth\'etie de rapport $\frac{1}{q}$
pour la distance $d_\infty$.

D'apr\`es \cite[Theo.~4.14]{Mat}, la dimension de Hausdorff $\delta$
de $(\O,d_\infty)$ v\'erifie l'\'equa\-tion $\sum_{\alpha\in k}\frac{1}{q}^\delta
=1$.  Donc $\delta =1$.

Comme par r\'ecurence, $\O$ est r\'eunion de $q^n$ boules (disjointes)
de rayon ${1\over q^n}$, la mesure de Hausdorff de $\O$ est au plus
$1$.  Rappellons que si deux boules pour $d_{\infty}$ se rencontrent,
alors l'une est contenue dans l'autre, et que toute boule contenue
dans $\O$ est une boule de rayon une puissance de ${1\over q}$. Donc
la mesure de Hausdorff de $\O$ vaut au moins $1$.
\cqfd

\bprop\label{prop:etape_deux}
Les mesures de Haar $\mu_{\rm Haar}$ et de Hausdorff $\mu_{\rm Haus}$
sur $\wh K$ co\"{\i}ncident.
\eprop

\dem Comme $\delta$ est finie non nulle, et $\O$ est un espace de
Cantor de mesure de Hausdorff finie, la mesure de Hausdorff $\mu_{\rm
     Haus}$ est une mesure de Radon (voir par exemple \cite[page
57]{Mat}).

Comme vu ci-dessus, la mesure $\mu_{\rm Haus}$ est invariante par
translation, donc par unicit\'e de la mesure de Haar, il existe $c\geq
0$ tel que $\mu_{\rm Haus}=c\mu_{\rm Haar}$.  Comme $\mu_{\rm
    Haar}(\O)=1$ et $\mu_{\rm Haus}(\O)=1$, les mesures co{\"\i}ncident.
\cqfd

\bigskip
Pour tous $\ga$ dans $G$ et $f$ dans $\wh K$, v\'erifiant
$f\neq\ga^{-1}\infty$,
notons
$$j_\infty(\ga,f)=|\ga'(f)|_\infty$$
la valeur absolue de la d\'eriv\'ee en $f$ de l'application holomorphe
$\ga:\wh K\ra\wh K$.  Par la formule de d\'erivation des compositions,
elle v\'erifie la formule de cocycle suivante~: pour $\ga_1,\ga_2$
dans $ G$ et $f$ dans $\wh K$ tels que $\ga_2f,\ga_1\ga_2f\neq \infty$,
$$j_\infty(\ga_1\ga_2,f)=j_\infty(\ga_1,\ga_2f)j_\infty(\ga_2,f)\;.$$
Si $\ga=\left[\begin{array}{cc} a & b \\ c & d\end{array}\right]$ avec
$|ad-bc|_\infty=1$, alors $j_\infty(\ga,f)=\frac{1}{|cf+d|^2_\infty}$.
Comme la valeur absolue prend des valeurs discr\`etes, pour tout
\'el\'ement $\ga$ de $G$, l'application continue $f\mapsto
j_\infty(\ga,f)$ est localement constante sur $\wh
K-\{\ga^{-1}\infty\}$.

\bprop\label{prop:etape_trois}
Pour tout $\ga$ de $G_1$, pour tous $f,g$ dans
$\wh{K}$ tels que $f,g\neq \gamma^{-1}\infty$,
$$|\ga f-\ga g|_\infty^2 =
j_\infty(\ga,f)j_\infty(\ga,g)|f - g|_\infty^2\;.$$
\eprop

\dem
Les translations $t_a=\left[\begin{array}{cc} 1 & a \\ 0 &
      1\end{array}\right]$ avec $a$ dans $\wh K$, les applications
$\lambda_{\alpha,\beta}=\left[\begin{array}{cc} \alpha & 0 \\ 0 &
      \beta \end{array}\right]$ avec $\alpha,\beta$ dans $\wh
K^{\times}$ tels que $|\beta|_\infty=|\alpha|_\infty^{-1}$, et
l'inversion $i=\left[\begin{array}{cc} 0 & 1 \\ 1 &
      0\end{array}\right]$ engendrent $G_1$. Par la formule de cocycle
pour $j_\infty$ et par un argument de continuit\'e sur les $f,g$, il
suffit donc de v\'erifier la proposition \ref{prop:etape_trois} pour
ces trois types de transformations.

Pour $a$ dans $\wh K$, la translation $t_a$ pr\'eserve la distance
$d_\infty$, et $j_\infty(t_a,f)=1$. Pour $\alpha,\beta$ dans $ \wh
K^{\times}$ avec $|\beta|_\infty=|\alpha|_\infty^{-1}$, l'application
$\lambda_{\alpha,\beta}$ est une homoth\'etie de rapport
$|\alpha|_\infty^{2}$ pour la distance $d_\infty$ et
$j_\infty(\lambda_{\alpha,\beta},f)=|\alpha|_\infty^2$. Comme
$j_\infty(i,f)=\frac{1}{|f|_\infty^2}$ et
$|\frac{1}{f}-\frac{1}{g}|_\infty= \frac{|f -
    g|_\infty}{|f|_\infty|g|_\infty}$, le r\'esultat en d\'ecoule.
\cqfd

\bprop\label{prop:etape_quatre}
Pour tout $\ga$ de $G_1$, l'image de la mesure $\mu_{\rm Haus}$ par
$\gamma$ est absolument continue par rapport \`a $\mu_{\rm Haus}$. Sa
d\'eriv\'ee de Radon-Nikodym v\'erifie (pour $\mu_{\rm Haus}$-presque
tout $\xi$ de $\wh K$):
$$\frac{d\ga_*\mu_{\rm Haus}}{d\mu_{\rm Haus}}(\xi)=
j_\infty(\ga^{-1},\xi)\;.$$
\eprop

\dem
Soit $\ga$ dans $ G$. Comme $f\mapsto j_\infty(\ga,f)$ est
localement constante, pour tout $f_0$ de $\wh K-\{\ga^{-1}\infty\}$,
il existe un voisinage $U$ de $f_0$ dans $\wh K-\{\ga^{-1}\infty\}$
tel que pour tout $f$ dans $ U$, on a
$j_\infty(\ga,f_0)=j_\infty(\ga,f)$. En particulier, la proposition
\ref{prop:etape_trois} montre que pour tous $f,g$ dans $ U$, on a
$|\ga f-\ga g|_\infty = j_\infty(\ga,f_0)|f - g|_\infty$.  Donc la
restriction de $\ga$ \`a $U$ est une homoth\'etie de rapport
$j_\infty(\ga,f_0)$ pour la distance $d_\infty$. Soit $V$ un voisinage
de $f_0$ d'adh\'erence contenue dans $U$. Par construction de la
mesure de Hausdorff, pour tout bor\'elien $E$ de $V$, on a $\mu_{\rm
    Haus}(\ga E)=j_\infty(\ga,f_0)^\delta\mu_{\rm Haus}(E)$.  Ceci
montre le r\'esultat car $\delta=1$.
\cqfd

\medskip
Nous aurons besoin du calcul \'el\'ementaire d'int\'egrale suivant.

\blemm \label{lem:calcintegral} Si $\mu=\mu_{\rm Haar}$, alors
$$\int_J \frac{d\mu(g)}{|g|^2_\infty}=\frac{1}{q}\;.$$
\elemm

\dem
Notons $S$ cette int\'egrale. Comme $J$ est la r\'eunion disjointe
des ensembles
$(a+X^{-1}\O)$, quand $a$ varie dans $A-k$, et comme la mesure de 
Haar est invariante par
translation,
$$S=\sum_{a\in A-k}\int_{X^{-1}\O}\frac{d\mu(g)}{|a+g|^2_\infty}
=\sum_{a\in A-k}\frac{1}{|a|^2_\infty}\mu(X^{-1}\O)\;.$$
Comme le nombre de polyn\^omes de degr\'e $n$ \`a coefficients dans
$k$ est $(q-1)q^n$,
$$S=\sum_{n\in\NN-\{0\}}\frac{(q-1)q^n}{q^{2n}}\frac{1}{q}=
\frac{1}{q}\sum_{n\in\NN-\{0\}}\frac{(q-1)}{q^{n}}=\frac{1}{q}\;.$$
\cqfd

\bigskip
\noindent{\bf D\'emonstration du th\'eor\`eme \ref{theo:mesinva}.}
Rappelons que 
$$d\wt m (\xi_-,\xi_+, x_0,v_0)=
\frac{d\mu(\xi_-)d\mu(\xi_+)dx_0dv_0}{|\xi_+ -\xi_-|^2_\infty}$$
et que les mesures de comptage sont invariantes par permutation. Il
est donc imm\'ediat que $\wt m$ est une mesure bor\'elienne positive,
invariante par renversement du temps $\wt\tau:(\xi_-,\xi_+,
x_0,v_0)\mapsto (\xi_+,\xi_-, x_0,v_0)$ et par $\wt T:(\xi_-,\xi_+,
x_0,v_0)\mapsto (\xi_-,\xi_+, x_1,v_1)$.

Montrons que $\wt m$ est invariante par $\Ga$. Toute action d'un
groupe sur un ensemble pr\'eserve la mesure de comptage. Pour tout
$\ga$ de $\Ga$, pour $\wt m$-presque tout $(\xi_-,\xi_+, x_0,v_0)$,
\begin{eqnarray*}
d(\ga_*\wt m) (\xi_-,\xi_+, x_0,v_0)
&=&
\frac{d(\ga_*\mu)(\xi_-)d(\ga_*\mu)(\xi_+)dx_0dv_0}{|\ga^{-1}\xi_+
     -\ga^{-1}\xi_-|^2_\infty}\\
&=&
\frac{j_\infty(\ga^{-1},\xi_-)d\mu(\xi_-)
     j_\infty(\ga^{-1},\xi_+)d\mu(\xi_+)dx_0dv_0}
{j_\infty(\ga^{-1},\xi_-)j_\infty(\ga^{-1},\xi_+)|\xi_+
     -\xi_-|^2_\infty}
\end{eqnarray*}
d'apr\`es les propositions \ref{prop:etape_quatre} et
\ref{prop:etape_trois}, et car $\Ga$ est contenu dans $G_1$.  Donc
$\ga_*\wt m=\wt m$.

\medskip
Par cons\'equent, $q^2\wt m$ induit par passage au quotient par $\Ga$
une mesure bor\'elienne positive $m$ sur $\G_0(\Ga\bac\TT_q)$,
invariante par $T$ et par renversement du temps $\tau$. Il ne reste
plus qu'\`a montrer que $m$ est une mesure de probabilit\'e.

Comme la mesure de Haar sur $\wh K$ est sans atome et comme $K$ est
d\'enombrable, le param\'etrage de $\G_0(\Ga\bac\TT_q)$ permet
d'\'ecrire~:
$$m(\G_0(\Ga\bac\TT_q))=q^2\int_{\xi_-\in J}\int_{\xi_+\in X^{-1}\O}
\frac{d\mu(\xi_-)d\mu(\xi_+)}{|\xi_+ -\xi_-|^2_\infty}=q^2\mu(X^{-1}\O)
\int_J \frac{d\mu(\xi_-)}{|\xi_-|^2_\infty}\;.$$
Le lemme \ref{lem:calcintegral} montre alors que la masse totale de
$m$ est bien
$1$.

Ceci ach\`eve la d\'emonstration du th\'eor\`eme \ref{theo:mesinva}.
\cqfd

\section{Applications}
\label{sec:applications}

Nous r\'esumons dans l'\'enonc\'e suivant les r\'esultats de la partie
\ref{sec:codage}. Toutes les applications en d\'ecouleront.

\bprop\label{prop:expropexprop_un}
Il existe un hom\'eomorphisme $\Theta'':\G_0(\Ga\bac\TT_q)\ra
(A-k)^{\ZZ}$ et une application continue surjective ${\rm
    p}_2:\G_0(\Ga\bac\TT_q)\ra X^{-1}\O\cap \,^c\! K$ tels que les
diagrammes suivants sont commutatifs~:
$$\begin{array}[b]{ccc} \G_0(\Ga\bac\TT_q) &
    \stackrel{\Theta''}{\longrightarrow} & (A-k)^{\ZZ}
    \\
    T\downarrow\;  & & \downarrow \sigma \\
    \G_0(\Ga\bac\TT_q) &\stackrel{\Theta''}{\longrightarrow} &
    (A-k)^{\ZZ}
\end{array}
\;\;\;\;\;\begin{array}[b]{ccc}
\G_0(\Ga\bac\TT_q) & \stackrel{T}{\longrightarrow} &
\G_0(\Ga\bac\TT_q)
\\
\downarrow {\rm p}_2 & & \downarrow {\rm p}_2 \\
X^{-1}\O\cap \,^c\! K &\stackrel{\Psi}{\longrightarrow} &
X^{-1}\O\cap \,^c\! K
\end{array}\;.$$
o\`u $\sigma$ est le d\'ecalage \`a gauche des suites bilat\`eres de
$(A-k)^{\ZZ}$.
\eprop

\dem Notons ${\rm p}_+:J\times X^{-1}\O\ra X^{-1}\O$ la seconde
projection. Avec les notations de la partie \ref{sec:codage}, posons
$$\Theta''=\Theta'\circ \varphi\;\;{\rm et}\;\;{\rm p}_2={\rm p}_+\circ
\varphi\;.$$
Le r\'esultat d\'ecoule alors de la partie
\ref{subsec:lienartin}.
\cqfd

\medskip
Le r\'esultat pr\'ec\'edent dit en particulier que
l'application d'Artin est semi-conju\-gu\'ee \`a l'application de
premier retour du flot g\'eod\'esique sur le rayon modulaire. C'est
l'analogue dans le corps des s\'eries de Laurent sur $\FF_q$ du
r\'esultat bien connu qui dit que l'application de Gauss $x\mapsto
\{\frac{1}{x}\}$ sur $[0,1]\cap\;^c\!\,\QQ$ est semi-conjugu\'ee \`a
une application de premier retour du flot g\'eod\'esique sur la courbe
modulaire ${\rm PSL}(2,\ZZ)\backslash\HH^2$ (voir par exemple
\cite{Seri}).

\subsection{Application au m\'elange du flot g\'eod\'esique
sur le rayon modulaire}
\label{sec:melange}

Consid\'erons la mesure $\nu$ sur l'espace discret $A-k$ d\'efinie par
$\nu(a)=\frac{1}{|a|_\infty^2}$. La mesure $\nu$ est une mesure de
probabilit\'e car le nombre de polyn\^omes de degr\'e $n$ \`a
coefficients dans $k$ est $(q-1)q^n$. Notons $\nu^\ZZ$ la mesure
produit sur $(A-k)^\ZZ$.

\bprop \label{prop:thetaisomesu}
L'hom\'eomorphisme $\Theta'':\G_0(\Ga\bac\TT_q)\ra (A-k)^\ZZ$
envoie la mesure $m$ sur la mesure $\nu^\ZZ$.
\eprop

\dem
Notons $P$ la mesure de probabilit\'e image de $m$ par $\Theta''$.
Comme $\Theta''$ conjugue $T$ et $\sigma$, et comme $m$ est invariante
par $T$, la mesure $P$ est donc invariante par le d\'ecalage $\sigma$.

Pour $a$ dans $ A-k$ fix\'e, calculons $P[a_{-1}=a]$. Soit $\ell=
(\xi_-,\xi_+) \in \G'_0(\TT_q)$.  Alors $a_{-1}(\ell)=a$ si et
seulement si $\xi_-$ appartient $ (a +X^{-1}\O)$ et $\xi_+$
appartient \`a $ X^{-1}\O$.  On
a donc~:
$$P[a_{-1}=a]= m((a +X^{-1}\O)\times X^{-1}\O)
=q^2\frac{1}{|a|^2_\infty}\mu(X^{-1}\O)^2=\nu(a)\;.$$

Soient $a'_{-k},\dots,a'_{-1},a'_{-0}$ fix\'es dans $A-k$. Par
s\'eparation des variables, les deux \'ev\`enements $\{a_{-k}=a'_{-k},
\dots, a_{-1}=a'_{-1}\}$ et $\{a_{0}=a'_{0}\}$ (qui portent
respectivement sur les coordonn\'ees $\xi_-$ et $\xi_+$) sont
ind\'ependants.

L'invariance de $P$ par le d\'ecalage $\sigma$, et un raisonnement par
r\'ecurrence montrent que les mesures $P$ et $\nu^\ZZ$ co{\"\i}ncident
sur tous les cylindres du produit $(A-k)^\ZZ$, donc sont \'egales.
\cqfd

\bcoro \label{coro:melangeflotgeod}
L'application $T:\G_0(\Ga\bac\TT_q)\ra \G_0(\Ga\bac\TT_q)$ de premier
retour du flot g\'e\-o\-d\'e\-si\-que sur le rayon modulaire est
Bernoulli (i.e.~(m\'etriquement) conjugu\'ee \`a un d\'ecalage de
Bernoulli sur un alphabet fini), donc m\'elangeante, donc ergodique,
pour la mesure $m$.
\ecoro

\dem
L'application $\Theta'':\G_0(\Ga\bac\TT_q)\ra (A-k)^\ZZ$ est un
isomorphisme mesur\'e, conjuguant $T$ au d\'ecalage $\sigma$. Le
d\'ecalage (de Bernoulli) $\sigma$ est d'entropie
$h_{\nu^\ZZ}(\sigma)$ finie, car
\begin{eqnarray*}
h_{\nu^\ZZ}(\sigma)
&=&
-\sum_{a\in A-k} \nu(a)\log\nu(a)= -2\log
q\sum_{a\in A-k} \frac{v_\infty(a)}{|a|_\infty^2}=2\log
q\sum_{n\in\NN}\frac{n(q-1)q^n}{q^{2n}}\\
&=&
2\log q(q-1) \frac{1}{q} \left( \frac{1}{1-\frac{1}{q}} \right)^2=
\frac{2q\log q}{q-1}\;.
\end{eqnarray*}

D'apr\`es \cite[Theo.~5, p.~53]{Orn}, un d\'ecalage (de Bernoulli)
sur un alphabet d\'enombrable, d'entropie finie, est (m\'etriquement)
conjugu\'e \`a un d\'ecalage (de Bernoulli) sur un alphabet fini.
Le r\'esultat en d\'ecoule.
\cqfd

\bcoro
Le temps moyen du $n$-\`eme parcours dans une horoboule de la
famille $\H\B$ d'une g\'eod\'esique d\'ecor\'ee de $\Ga\bac\TT_q$
est $T_n=\frac{2q}{q-1}$.
\ecoro

\dem Dans \cite{Pau},  il est montr\'e qu'une g\'eod\'esique d\'ecor\'ee
rencontre sa $n$-\`eme horoboule de la famille $\H\B$ sur un segment
de longueur $-2v_\infty(a_n)$. L'invariance par le d\'ecalage montre que
$$
T_n = T_0= E[-2v_\infty(a_0)]=-2\sum_{a\in A-k}
v_\infty(a)P[a_0=a]= -2\sum_{a\in A-k}
\frac{v_\infty(a_0)}{|a|^2_\infty}=\frac{2q}{q-1}\;.
$$
\cqfd

\subsection{Applications arithm\'etiques}
\label{sec:appliarithm}

Le r\'esultat suivant est l'analogue dans le corps des s\'eries
formelles de Laurent du r\'esultat bien connu qui dit que la mesure de
Gauss $\frac{1}{\log 2}(\frac{dt}{1+t})$ sur $[0,1]$ est la projection
de la mesure de Liouville sur le fibr\'e unitaire tangent \`a la
courbe modulaire ${\rm PSL}(2,\ZZ)\backslash\HH^2$ (voir par exemple
\cite{Seri}).

\bprop\label{prop:exprop_deux} L'image de $m$ par $p_2$ est la
restriction de la mesure $q\mu_{\rm Haar}$ \`a $X^{-1}\O$~:
$$(p_2)_*m=q \;\mu_{\rm Haar}\mid_{X^{-1}\O\cap\,^c\! K}\;.$$
\eprop

\dem
Par d\'efinition de $m$, pour tout $ \xi_+$ dans
$X^{-1}\O\cap\,^c\!  K$, nous avons
$$d(p_2)_*m(\xi_+)=q^2\int_{\xi_-\in J}\frac{d\mu(\xi_-)}
{|\xi_+-\xi_-|^2_\infty} d\mu(\xi_+)=q\,d\mu(\xi_+)$$
d'apr\`es le lemme \ref{lem:calcintegral}. D'o\`u le r\'esultat.
\cqfd

\medskip
Nous retrouvons ainsi l'invariance bien connue de la mesure
de Haar par l'application d'Artin, voir par exemple \cite{BN} pour des
r\'ef\'erences.

\bcoro\label{invarHaarparArtin} La mesure de Haar sur
$X^{-1}\O\cap\,^c\! K$ est invariante par l'application d'Artin.
\ecoro

\dem Les propositions \ref{prop:expropexprop_un} et
\ref{prop:exprop_deux}, et l'invariance de $m$ par l'application
de premier retour $T$ montrent que~:
$$\Psi_*\mu_{\rm Haar}=\frac{1}{q}(\Psi\circ p_2)_*m=\frac{1}{q}(
p_2\circ T)_*m = \frac{1}{q} (p_2)_*m = \mu_{\rm Haar}\;.$$
Le r\'esultat en d\'ecoule.  \cqfd

\medskip Nous renvoyons par exemple \`a \cite{BN} pour d'autres
propri\'et\'es dynamiques du d\'e\-ve\-lop\-pe\-ment en fractions
continues dans le corps des s\'eries de Laurent sur $\FF_q$.

\bigskip
\noindent {\small
\begin{tabular}{l}
Laboratoire de Math\'ematique UMR 8628 CNRS\\
Equipe de Topologie et Dynamique (B\^at. 425)\\
Universit\'e Paris-Sud \\
91405 ORSAY Cedex, FRANCE.\\
{\it e-mail: Anne.Broise@math.u-psud.fr}
\end{tabular}
\\
    \mbox{}
\\
    \mbox{}
\\
\begin{tabular}{l}
D\'epartement de Math\'ematique et Applications, UMR 8553 CNRS\\
Ecole Normale Sup\'erieure\\
45 rue d'Ulm\\
75230 PARIS Cedex 05, FRANCE\\
{\it e-mail: Frederic.Paulin@ens.fr}
\end{tabular}
}

\end{document}